\documentclass{article}

\usepackage{arxiv}

\usepackage[utf8]{inputenc} 
\usepackage[T1]{fontenc}    
\usepackage{hyperref}       
\usepackage{url}            
\usepackage{booktabs}       
\usepackage{amsfonts}       
\usepackage{nicefrac}       
\usepackage{microtype}      
\usepackage{lipsum}		
\usepackage{graphicx}
\usepackage{algorithm, algorithmic}
\usepackage{subcaption}
\usepackage{adjustbox}
\usepackage{doi}
\usepackage{amsmath}
\usepackage{xcolor}
\usepackage{xparse}
\usepackage{multirow}
\usepackage{tikz, array}

\usetikzlibrary{mindmap,trees,arrows,shapes,backgrounds,matrix,decorations.pathreplacing,decorations.pathmorphing,positioning}
\tikzstyle{every picture}+=[remember picture]

\everymath{\displaystyle}
\tikzstyle{na} = [baseline=-.5ex]

\DeclareDocumentCommand{\expect}{o o}{%
    \mathbb{E}%
    \IfNoValueTF{#2}%
        {%
            \IfNoValueTF{#1}%
            {}%
            { \left [ #1 \right ]}%
        }%
        {_{#1} \left [ #2 \right]}
}

\DeclareDocumentCommand{\probMeas}{o o}{%
    \probSym%
    \IfNoValueTF{#2}%
        {%
            \IfNoValueTF{#1}%
            {}%
            { \left ( #1 \right)}%
        }%
        {_\text{{#1}} \left (#2 \right )}
}

\title{An Autoencoder Compression Approach for Accelerating Large-scale Inverse Problems}

\ifx\anonymousAuthor\undefined
\author{ \hspace{1mm}Jonathan ~Wittmer \\
	Oden Institute for Computational Science and Engineering\\
	University of Texas at Austin\\
	Austin, TX 78712 \\
	\texttt{jonathan.wittmer@utexas.edu} \\
	\And
	\hspace{1mm}Jacob ~Badger \\
	Oden Institute for Computational Science and Engineering\\
	University of Texas at Austin\\
	Austin, TX 78712 \\
	\texttt{jcbadger@utexas.edu} \\
	\And
	\hspace{1mm}Hari ~Sundar \\
	Department of Computer Science\\
	University of Utah\\
	Salt Lake City, Utah 84112\\
	\texttt{hari.sundar@utah.edu}
	\And
	\hspace{1mm}Tan ~Bui-Thanh \\
	Department of Aerospace Engineering and Engineering Mechanics\\
	Oden Institute for Computational Science and Engineering\\
    University of Texas at Austin\\
	Austin, TX 78712 \\
	\texttt{tanbui@oden.utexas.edu}
}
\fi

\renewcommand{\shorttitle}{Autoencoder compression for inverse problems}
\hypersetup{
pdftitle={autoencoder_compression_inverse_problems},
pdfauthor={Jonathan ~Wittmer, Jacob ~Badger, Hari ~Sundar, Tan ~Bui-Thanh},
pdfkeywords={Machine Learning, Autoencoder, Inverse Problems, High Performance Computing},
}

\newcommand{\norm}[1]{\left\| #1 \right\|}

\newcommand{\LRp}[1]{\left( #1 \right)}
\newcommand{\LRs}[1]{\left[ #1 \right]}

\newcommand{\LRc}[1]{\left\{ #1 \right\}}

\newcommand{\bs}[1]{\boldsymbol{#1}}

\newcommand{\mb}[1]{\bs{#1}}

\newcommand{\bigO}{\mathcal{O}}

\newcommand{\half}{\frac{1}{2}}

\newcommand{\grad} {\ensuremath{\nabla}}
\newcommand{\Div} {\ensuremath{\nabla}\cdot}

\newcommand{\xb}{\bs{x}}
\newcommand{\R}{\mathbb{R}}

\newcommand{\y}{{y}}
\newcommand{\yb}{{\mb{\y}}}

\newcommand{\eval}[2][\right]{\relax
  \ifx#1\right\relax \left.\fi#2#1\rvert}

\newcommand{\Cb}{{\bs{C}}}

\newcommand{\Wb}{{\bs{W}}}

\newcommand{\argmin}{\text{arg min }}

\newcommand{\strain}{e}
\newcommand{\velocity}{\bs{v}}
\newcommand{\density}{\rho}
\newcommand{\cp}{c}
\newcommand{\partialD}[2]{\frac{\partial #1}{\partial #2}}
\newcommand{\rhs}{\bs{g}}
\newcommand{\domain}{\Omega}
\newcommand{\boundary}{\Gamma}
\newcommand{\inv}[1]{{#1}^{-1}}

\newcommand{\qoi}{\bs{u}}
\newcommand{\pto}{\bs{\mathcal{F}}}
\newcommand{\noiseCov}{\Gamma_{\text{noise}}}
\newcommand{\invNoiseCov}{\inv{\Gamma}_{\text{noise}}}
\newcommand{\priorCov}{\Gamma_{\text{prior}}}
\newcommand{\invPriorCov}{\inv{\Gamma}_{\text{prior}}}

\newcommand{\obs}{\bs{d}}

\newcommand{\ident}{\mathcal{I}}

\newcommand{\Wbo}{\bs{W}_1}
\newcommand{\Wbt}{\bs{W}_2}
\newcommand{\proj}{\bs{P}_2}
\newcommand{\qoiDias}{\qoi_{\text{DIAS}}}

\newcommand{\probSym}{\pi}
\newcommand{\loss}{\mathcal{L}}
\newcommand{\given}{\vert}

\begin{document}
\bibliographystyle{apalike}
\maketitle

\begin{abstract}
    PDE-constrained inverse problems are some of the most challenging and computationally demanding 
problems in computational science today. Fine meshes that are required to accurately
compute the PDE solution introduce an enormous number of parameters and require large
scale computing resources such as more processors and more memory to solve such systems in a reasonable time.
For inverse problems constrained by time dependent PDEs, the adjoint method 
that is often employed to efficiently compute gradients and higher order derivatives
requires solving a time-reversed, so-called adjoint PDE that depends on the forward
PDE solution at each timestep. This necessitates the storage of a high dimensional forward solution
vector at every timestep. 
Such a procedure quickly exhausts the available memory resources.
Several approaches that trade additional computation for reduced memory footprint
have been proposed to mitigate the memory bottleneck, including checkpointing
and compression strategies. In this work, we propose a 
close-to-ideal scalable compression approach using autoencoders to eliminate the need for 
checkpointing and substantial memory storage, thereby reducing both the time-to-solution and  memory requirements.
We compare our approach with checkpointing and  an off-the-shelf compression approach
on an earth-scale ill-posed seismic inverse problem. The results verify the expected close-to-ideal speedup for both the gradient and Hessian-vector product using the proposed autoencoder compression approach. 
To highlight the usefulness of the proposed approach, we combine the  autoencoder compression 
with the data-informed active subspace (DIAS) prior to show how the DIAS method
can be affordably extended to large scale problems without the need of checkpointing and large memory.
 
\end{abstract}


\graphicspath{{original_content/}}

\section{Introduction}
Modern supercomputers are an essential tool for today's scientists and engineers, enabling 
them to simulate larger and more complex systems than ever before. The availability of 
large scale computing resources has enabled numerous breakthroughs in scientific knowledge 
in areas such as quantum computing \cite{liu2021closing,doi2019quantum,mandra2021hybridq},
computational chemistry \cite{de2010utilizing,kowalski2021nwchem},
drug discovery \cite{ge2013molecular,schmidt2017next,bharadwaj2021computational,sukumar2021convergence}, 
biology \cite{mcfarlane2008beatbox,bukowski2010biohpc},
and plasma physics \cite{bhattacharjee2021preface,fedeli2022pushing}.
While the increased number of floating-point operations per second (FLOPS) has enabled larger and 
more challenging problems to be solved in reasonable amounts of time, 
applications are increasingly becoming bottlenecked by limited memory of computing systems
--- especially those applications dealing with \textit{big data}
\cite{denis2022modeling,imani2019digitalpim}. 

There are at least three common approaches to mitigating this memory bottleneck:
develop new algorithms that require less storage (for example, using density functional theory 
to approximate electron interactions rather than directly solving a many body problem in 
quantum mechanics \cite{bickelhaupt2000kohn}), 
utilize additional storage 
devices such as hard drives to increase the amount of data that can be stored,
and trade extra computation for reduced storage requirements. 
Algorithmic improvement is always desirable, but is not always possible. Additionally, 
solving a closely related problem that has computational or storage advantages 
may cause a loss in accuracy. 
Serializing data to disk was previously only a viable solution for problems with 
high compute intensity.
The high cost of reading from and writing to disk can be hidden when sufficient
computation is performed. 
However, recent 
advances in storage technology such as non-volatile memory
(NVM) have made serialization approaches more 
viable \cite{peng2020demystifying}. Indeed, NVM powers the large-memory nodes on 
Frontera, the Texas Advanced Computing Center's largest supercomputer at the time of writing.
This allows for traditional DRAM to be used as an
extra level of cache \cite{wu2017early,stanzione2020frontera}.
The last mentioned approach, trading computation for storage, 
includes checkpointing methods \cite{wang2009minimal,zhang2023optimal,griewank1997treeverse}
and sparse decomposition \cite{zhao2020smartexchange}.
More generally, 
this means storing the inputs required to regenerate the desired output rather than storing 
the output directly \cite{akturk2018trading}. In addition to saving time,
it is shown in \cite{akturk2018trading}
that recomputation can have energy saving benefits over storing and retrieving data. 

One type of problem that calls for advanced big data management techniques is 
full-waveform seismic inversion (FWI). As motivation, we give a high-level 
description of why advanced data management techniques are required for FWI 
here and give a detailed description in Section \ref{sec:problem_description}.
FWI is a partial differential equation (PDE) constrained inverse problem 
with dependencies on both space and time derivatives. The dependency on 
time leads to the generation of large amounts of data. With spatio-temporal 
measurements of the velocity field, the inverse problem is
to infer the underlying spatially varying material properties, namely the 
acoustic wave speed. Typical inverse solution methods reframe the problem 
as a constrained optimization problem and use gradient-based methods 
to minimize an objective function.
Employing the adjoint approach \cite{giles2000introduction,fichtner2006adjoint}
to compute the gradient necessitates solving an adjoint PDE that 
is solved backward in time and depends on the corresponding forward 
solution at each timestep. As a result, \textit{the entire 
solution history of the forward problem is needed in reverse order 
to solve the adjoint equation}. For large scale problems with high
spatial resolution and many timesteps, storing the entire forward solution 
history is infeasible, requiring a different approach. 
As a quick back-of-the-envelope calculation, 
consider a domain with 1 billion degrees of freedom and 2 000 timesteps. Using a 4$^{th}$ order 
Runge-Kutta scheme to integrate in time, there are 4 solutions to be stored per timestep. Each 
solution has 6 fields --- 3 velocity fields and 3 strain fields
(assuming a velocity-strain formulation).
This would require 384 TB to store in double precision.
Perhaps one of the 
most commonly used techniques for such a problem is to employ checkpointing. 
Checkpointing is a technique whereby solutions are stored at only 
a select few timesteps and the rest discarded.
The solution at any timestep can then be recreated by solving the forward 
problem from the nearest previous checkpoint.
While reducing memory requirements, checkpointing effectively increases the
computation requirements by 1 PDE solve per gradient evaluation. 

An alternative approach employs compression to reduce both the memory requirements compared to full storage
and computation requirements compared to the checkpointing case.
The idea to use compression to mitigate the effects of limited memory or storage is not new. 
Compression techniques have been used extensively in the audio and video processing community 
for decades \cite{blesser1969audio,lewis1990video}.
There are two kinds of compression: lossless compression
and lossy compression. Lossless compression is a class of techniques that are able to 
exactly recover the original input data. The PNG image format is an example of a format 
that employs lossless compression \cite{kaur2016review}.
The checkpointing strategy can be viewed as a form 
of lossless compression where information is encoded in the 
checkpoints and the physics model provides a decompression algorithm.
Lossy compression, on the other hand, allows for errors to be made in reconstructing 
the input data. 
The JPEG protocol is an example of lossy compression
in the image processing domain \cite{usevitch2001tutorial}. 

Various compression techniques have also been proposed to aid 
in scientific computation. Specific to the seismic inversion community, there are two 
main avenues of compression: compressing the source/receiver data \cite{habashy2011source} 
and compressing the forward solution to avoid checkpointing. Source/receiver compression
is useful in the case where large amounts of data are available, such as when many shots are 
recorded. The gradient must be evaluated at each shot;
each requires the solution of the forward and 
adjoint PDEs, implying that the cost scales linearly in the number of shots \cite{duarte2020seismic}.
The goal of source/receiver compression is to find linear combinations
of simultaneous sources and receivers
that are (nearly) equivalent to the original problem. When the number of linear combinations required 
to closely replicate the original problem is small, the number of forward and adjoint solves required 
to compute a single gradient is much smaller, reducing both the memory footprint and
the computation time \cite{habashy2011source}. 

Several other works have already explored the idea of using compression
to mitigate high adjoint-induced memory requirements
\cite{cyr2015towards,kukreja2019combining,kukreja2022lossy,boehm2016wavefield}.
Those works demonstrate that data can be compressed
without negatively affecting the inverse solution and
that compression can be faster than forward simulation.
Several techniques of compressing in space 
as well as a spline interpolation technique for compression in time
were proposed in \cite{boehm2016wavefield}.
Results for a variety of
off-the-shelf compression algorithms, including extensive results
using the ZFP compression algorithm \cite{lindstrom2014fixed},
were given in \cite{kukreja2019combining}.
Each of these papers demonstrates how different compression
techniques can be used to solve seismic inverse problems on
simple domains with large amounts of data.
That is, the inverse problem is not 
ill-posed and regularization techniques are not required.
While perturbations to the gradient due to compression may
minimally affect the inverse solution in the 
case of large data, it is unclear that such approaches will be viable
in solving ill-posed problems with limited data. Additionally,
these approaches have limited compression capabilities,
limiting the maximum problem size that can be solved. 

\begin{figure}[!htp]
  \centering
  \begin{tikzpicture}[
    plain/.style={
      draw=none,
      fill=none,
    },
    net/.style={
      matrix of nodes,
      nodes={
        draw,
        circle,
        inner sep=8.5pt
      },
      nodes in empty cells,
      column sep=0.2cm,
      row sep=-11pt
    },
    >=latex
    ]
    \matrix[net] (mat)
    {
      |[plain]| \parbox{1cm}{\centering Input\\layer} & |[plain]| & |[plain]| \parbox{1cm}{\centering Latent\\space}
      & |[plain]| & |[plain]| \parbox{1.5cm}{\centering Output\\layer} \\
      & |[plain]| & |[plain]| & |[plain]| & \\
      |[plain]| & |[plain]| & |[plain]| & |[plain]| & |[plain]|\\
      & & |[plain]| & & \\
      |[plain]| & |[plain]| & & |[plain]| & |[plain]| \\
      & & |[plain]| & & \\
      |[plain]| & |[plain]| & & |[plain]| & |[plain]| \\
      & & |[plain]| & & \\
      |[plain]| & |[plain]| & |[plain]| & |[plain]| & |[plain]|\\
      & |[plain]| & |[plain]| & |[plain]| & \\
    };

    \foreach \ai in {2,4,6,8,10}
    {\foreach \aii in {4,6,8}
      \draw[->] (mat-\ai-1) -- (mat-\aii-2);
    }

    \foreach \ai in {4,6,8}
    {\foreach \aii in {5,7}
      \draw[->] (mat-\ai-2) -- (mat-\aii-3);
    }

    \foreach \ai in {5,7}
    {\foreach \aii in {4,6,8}
      \draw[->] (mat-\ai-3) -- (mat-\aii-4);
    }

    \foreach \ai in {4,6,8}
    {\foreach \aii in {2,4,6,8,10}
      \draw[->] (mat-\ai-4) -- (mat-\aii-5);
    }

    \draw[
        thick,
        decoration={
          brace,
          mirror,
          raise=0.5cm,
          amplitude=10pt
        },
        decorate,
        anchor=south,
    ] (mat-10-1) -- (mat-10-3)
    node [pos=0.5,anchor=north, yshift=-0.75cm] {Encoder};

    \draw[
        thick,
        decoration={
          brace,
          mirror,
          raise=0.5cm,
          amplitude=10pt
        },
        decorate,
        anchor=south,
    ] (mat-10-3) -- (mat-10-5)
    node [pos=0.5,anchor=north, yshift=-0.75cm] {Decoder};
    
  \end{tikzpicture}
  
  \caption[Simple autoencoder]{Autoencoder architecture with deep neural networks}
  \label{fig:autoencoder_simple}
\end{figure}
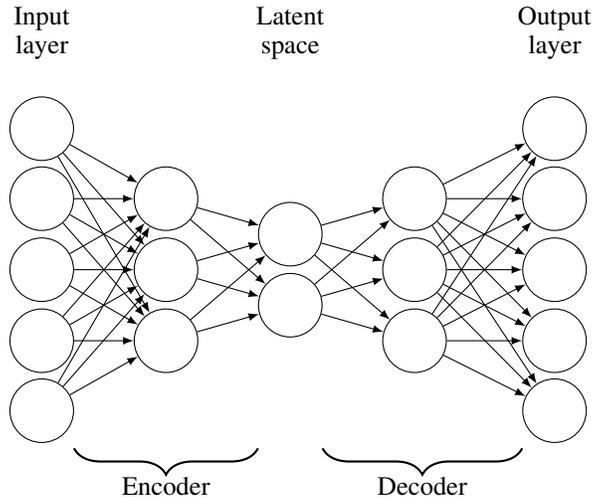

Meanwhile, the scientific and technology communities have seen
an explosion in the applications of machine learning techniques.
Deep learning in particular has found fruitful application in enabling
progress in speech recognition \cite{kamath2019deep},
computer vision \cite{bjerge2022real,esteva2021deep},
and in accelerating innovation in scientific computation
\cite{reichstein2019deep,kates2019predicting}.
One approach that has emerged is to use an autoencoder architecture
to perform data compression \cite{cheng2018deep,liu2021high,hinton2006reducing}.
An autoencoder is simply a deep neural network
that approximates the identity map.
The output at some intermediate layer can be extracted as a
representation of the input and is referred to 
as the \textit{latent representation} which 
resides in the \textit{latent space} \cite{dillon2021better}.
Autoencoders are typically decomposed into 
two separate networks: an encoder which consists of the layers
from the input to the latent space and a decoder which is a network
that maps the latent space back to the original dimension.
To achieve a compressed representation, the latent space must have dimension strictly smaller
than the input. Compared to traditional dimension reduction methods,
autoencoders can be viewed as a non-linear version of principal component analysis (PCA) \cite{ladjal2019pca,kneer2021symmetry}.
In the most simple case of a single layer encoder and single layer decoder with linear activation,
the autoencoder learns the same compression as PCA. That is, the single layer linear autoencoder projects the
input data onto the same lower dimensional subspace that PCA identifies,
though the representation is different \cite{baldi1989neural,ladjal2019pca}.
The actual PCA basis vectors can be obtained from the weights of the trained encoder,
showing the two are equivalent \cite{plaut2018principal}.
The advantage of a deep autoencoder with nonlinear activations is that it is often able to achieve lower reconstruction error for
the same dimension latent space compared to PCA \cite{liu2021exploring,liu2021high},
enabling higher feasible compression ratios in application.

While autoencoders have empirically demonstrated impressive compression ratios and low error rates,
there is no theory that we are aware of at this time bounding the error of an autoencoder.
However, a practical method of bounding the error of an autoencoder is to compute and store some form of 
the residual between the desired output ($\bs{y}_{\text{true}}$) and the actual output ($\bs{y}_{\text{pred}}$),
i.e. $\bs{r} = \bs{y}_{\text{true}} - \bs{y}_{\text{pred}}$. Various techniques of approximately representing
this residual vector have been proposed such as compressing $\bs{r}$ with an off-the-shelf
algorithm when $\norm{\bs{r}} > \text{tol}$ for some tolerance, $\text{tol}$ \cite{lee2022error,liu2021exploring} or
storing only the elements $r_i$ of $\bs{r}$ where $|r_i| > \text{tol}$ \cite{liu2021high,abu2014efficient}.

To address the problems of additional expensive PDE solves required by the checkpointing approach
and the limited compression capability of off-the-shelf compression techniques,
we propose a novel autoencoder approach for compressing the forward solution.
We show empirically that this compression method is faster than the checkpointing 
method while achieving comparable solution accuracy.
In addition to showing that a trained autoencoder has memory and computational advantages
over existing methods, we develop an efficient data generation and training procedure
based on the Bayesian inversion formulation that enables our approach
to maintain its advantages, even when including the cost of training the autoencoder. 
While variations of autoencoders are often used
for generative modeling \cite{kingma2013auto}
and others assign a physical interpretation
to the latent space \cite{goh2019solving},
we use the \textit{auto-encoding}
capabilities of the autoencoder architecture to allow the training process
to find the most effective compressed latent representation for seismic input data.
We provide a mathematical description of the FWI problem and the inverse problem in
Section \ref{sec:problem_description}. 
Two autoencoder compression variants are proposed in Section \ref{section:ae_methodology}
along with detailed explanation of autoencoder architecture, training data selection, and normalization.
We show in Section \ref{section:ae_numerical_results} that machine learning can be leveraged
to accelerate the solution of seismic inverse problems
in very high dimensions, even in the ill-posed case with few measurements.
Through numerical experiment, 
we show that this approach is faster than the checkpointing approach, and as fast as
the state-of-the-art off the shelf compression approach, ZFP \cite{lindstrom2014fixed},
while achieving much higher compression ratios.
We then extend the data informed active subspace (DIAS) regularization method 
\cite{nguyen2022dias}
to nonlinear problems and show that autoencoder compression can be used to enable the
affordable application of the DIAS prior for large-scale inverse problems.


\section{Full Waveform Inversion}
\label{sec:problem_description}
As a simplified model of the earth's seismic properties, we consider the
acoustic wave equation. In velocity-strain form, this gives rise to a
system of first order PDEs. We choose the velocity-strain form because
of its close relation to the acoustic-elastic wave equation \cite{wilcox2015discretely}.
This enables simultaneous derivation of gradients and Hessians
for both acoustic, elastic, and acoustic-elastic formulations.
Additional computational and mathematical advantages can
be found in \cite{wilcox2015discretely}.
Following the setup given in \cite{bui2013computational}, the acoustic
wave equation is given by, 
\begin{subequations}
  \label{eqn:model}
  \begin{align}
    \density \partialD{\velocity}{t} - \grad \LRp{\density \cp^2 \strain} = \rhs, \\
    \partialD{\strain}{t} - \Div \velocity = 0,
    \label{eqn:dilatation_time_derivative}
  \end{align}
\end{subequations}
where $\density = \density(\xb)$ is the density, $\cp = \cp(\xb)$ is the acoustic
wave speed for which we are inverting, $\rhs = \rhs(\xb, t)$ is a forcing function,
$\velocity = \velocity(\xb, t)$ is the velocity vector, $\strain = \strain(\xb, t)$
is the trace of the strain tensor (dilatation),
and $\xb$ denotes the spatial coordinate in the domain $\domain$.
We specify initial conditions for the velocity and dilatation fields
\begin{equation}
  \label{eqn:initial_conditions}
  \velocity(\xb, 0) = \velocity_{0}(\xb) \quad \text{and} \quad
  \strain(\xb, 0) = \strain_{0}(\xb),  \quad \xb \in \domain.
\end{equation}
Lastly, we impose traction-free boundary conditions
\begin{equation*}
  \label{eqn:boundary_conditions}
  \strain(\xb, t) = 0, \quad \xb \in \boundary = \partial \domain, \; t \in (0, T)
\end{equation*}
with final time $T$.
The 3 components of the velocity vector field and the 3 diagonal components of
the strain tensor field, used to compute the dilatation, will be collectively
referred to as the \textit{state} variables in the rest of this paper.
It is the state variables that need to be compressed and decompressed. 

In this work, we consider the inverse parameter estimation problem ---
inferring the wave speed $\cp\LRp{\xb}$
from sparse measurements of the velocity, $\obs$.
Since we consider the setting with limited measurements, the inverse problem is
ill-posed and regularization is required. The inverse problem can be formulated as the
optimization problem
\begin{equation}
  \qoi^* := \underset{\qoi}{\argmin} \half \norm{\pto \LRp{\qoi} - \obs}_{\invNoiseCov}^2
  + \half \norm{\qoi - \qoi_0}_{\invPriorCov}^2
\end{equation}
where $\qoi$ is referred to as the parameter of interest (PoI),
$\pto$ is the parameter to observable map (PtO map), $\noiseCov$ is the
noise covariance matrix of the observations $\obs$, $\qoi_0$ is the initial guess
and $\priorCov$ characterizes the regularization. In the Bayesian setting
\cite{kaipio2006statistical},
$\qoi_0$ and $\priorCov$ are interpreted as the prior mean and prior covariance matrix, respectively.
While we consider full waveform inversion in the deterministic setting,
it is convenient to consider the problem in the Bayesian framework and
only work toward estimating the maximum a posterior (MAP) point.
The statistical framing of the MAP problem allows us to choose regularization
via selection of a prior distribution rather than ad-hoc and allows
for trivial extension of our proposed compression approach to quantifying uncertainty. 
In our case, the PoI is the acoustic wave speed which is mapped to our observations, $\obs$, 
by solving the acoustic wave equation and observing the velocity field only at the receivers.
Let $\qoi \approx \cp$ be the discretized acoustic wave speed that we will numerically estimate.
Then
\begin{equation*}
  \pto := \bs{B} \bs{\mathcal{A}}
\end{equation*}
where $\bs{\mathcal{A}}$ is the solution of \eqref{eqn:model} and $\bs{B}$ is an
observation operator, extracting the velocity field at receiver locations. 

As in \cite{bui2013computational}, we choose the prior to be the PDE-based BiLaplacian
which encodes smoothness and anisotropy of the wave speed along with our certainty about our initial guess.
The BiLaplacian prior term $\half \norm{\qoi - \qoi_0}_{\invPriorCov}^2$
is computed by solving the following elliptic PDE:
\begin{subequations}
\begin{align}
  -\alpha \grad \cdot \LRp{\Theta \grad \qoi} + \alpha \qoi = s \quad &\text{in } \Omega\\
  \alpha \LRp{\Theta \grad \qoi} \cdot \bs{n} = 0 \quad &\text{on } \partial \Omega
\end{align}
\end{subequations}
where $\bs{n}$ is the outward facing unit normal on the boundary, $\partial \Omega$.  
$\invNoiseCov$ is the square of this differential operator. Properties of this
prior and further discussion can be found in \cite{bui2013computational}.
Finally, derivation of the discrete forward, adjoint, and gradient expressions
for a discontinuous Galerkin discretization can be
found in \cite{wilcox2015discretely}.

\section{Compression methods}
\label{section:ae_methodology}
There are a variety of ways that the state variables
can be compressed. Each state variable is a function of time and
3 space dimensions, resulting in 4 dimensions across which
the state can be compressed. Additionally, there might be some relationship
between the state variables that could be exploited to achieve higher
compression ratios. In this work, we propose two 
compression techniques that are amenable for implementation with
autoencoders, treating each state variable independently: 
\begin{enumerate}
\item Compression in space, allowing the autoencoder to discover the 3D structure of the state across multiple elements.
\item Compression in time on an element-by-element basis.
\end{enumerate}

As a brief justification for pursuing only these two compression setups, let us
make a few notes. First, there is potential for scale mismatch between
each of the state variables. Consider again equation \eqref{eqn:dilatation_time_derivative}.
The time derivative of the dilatation is related to the divergence of the
velocity field. From this we can see that there could be very large velocities,
but no dilatation if the velocity field is divergence free.
This makes it difficult to design a method that can garner any
correlation between the \textit{value} of the velocity field at a given timestep
and the \textit{value} of the strains.

Secondly, we decided against pursuing methods requiring explicit specification
of the fields in three-dimensional coordinates. That is, we treat all the
nodal coefficients equally, not explicitly taking advantage of the embedding
of these nodal coefficients in 3-dimensional space. While it is in principle
possible to use something like a convolutional neural network with
3-dimensional kernels \cite{ji20123d}, we found two major difficulties:
1) convolutional neural networks in 3D are rather slow and
2) the nodes of our DG mesh are not equally spaced.
Although equally spaced nodes would enable simple application of
convolutional autoencoders, non-uniformly spaced nodal coordinates
such as Gauss-Lobatto-Legendre (GLL) enable more accurate integration \cite{quarteroni2010numerical}
and are a natural choice for node location in a DG scheme. 
In Sections \ref{section:space_compression} and \ref{section:time_compression},
we will present the implementations of the two proposed autoencoder compression approaches. 

\subsection{Compression on a uniform mesh}
\label{section:space_compression}
\begin{figure}[htb!]
  \centering
  \begin{subfigure}{0.4\textwidth}
    \includegraphics[trim=1cm 1cm 4.5cm 1cm,clip=true,width=\textwidth]{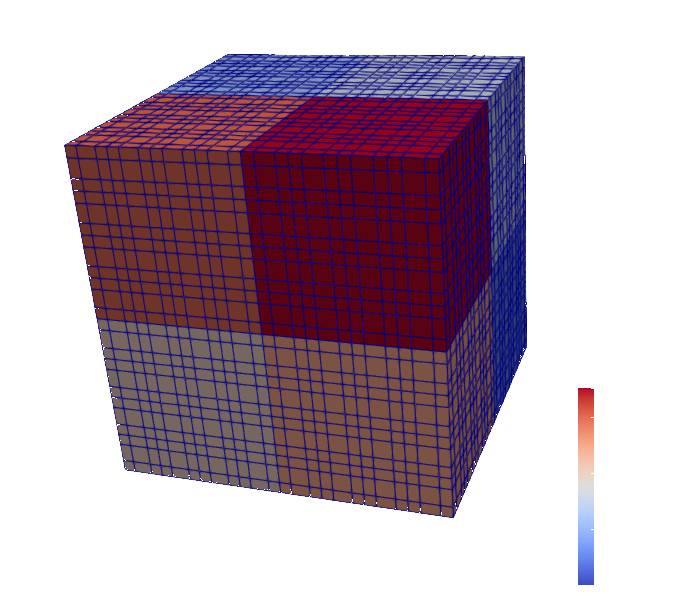}
    \caption{}
  \end{subfigure}
  \hspace*{1cm}
  \begin{subfigure}{0.4\textwidth}
    \includegraphics[trim=1cm 1cm 4.5cm 1cm,clip=true,width=\textwidth]{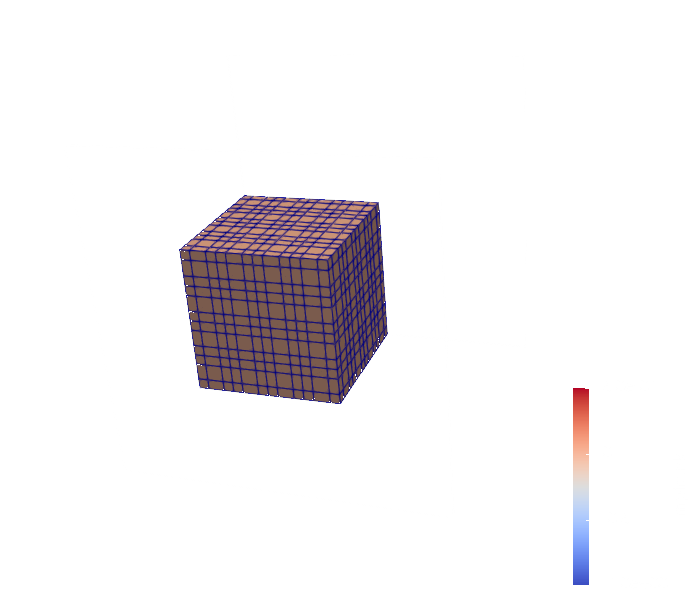}
    \caption{}
  \end{subfigure}
  \caption[Box domain mesh]{
    (a) Course mesh on box domain showing non-uniform spacing of grid.
    Different colors correspond to different MPI ranks.
    (b) mesh of single MPI rank consisting of 8 third order elements. 
  }
  \label{fig:box_mesh}
\end{figure}

Consider a box domain where $\domain = \LRs{0, 1}^3$.
Such a domain may arise in reservoir exploration,
the example used in \cite{kukreja2019combining}.
Figure \ref{fig:box_mesh} shows the mesh on a box domain.
While the nodes are unevenly spaced, they occur in a regular pattern.
This allows us to flatten a state variable across multiple elements
into a single vector for compression.
Care must be taken to partition the domain evenly across MPI ranks and to ensure that 
each rank stores the nodes in the same spatial orientation, i.e.
the elements are stacked the same way in memory on each rank. 
The autoencoder then has the freedom to identify the most relevant spatial
correlation in order to most accurately compress the data.
Note that the order in which the elements are stacked does not matter
so long as it is the same across all ranks since a permutation can be implicitly learned. 
With the same state layout on every rank,
a single autoencoder can be trained that compresses the
local degrees of freedom on a single rank. Once trained,
each rank stores its own copy of the autoencoder, eliminating the need for any communication.
The fact that no communication is required is key to the scalability of our proposed approach. 

There are a few competing factors that must be balanced in designing an autoencoder for compression.
The input vector should be:
\begin{enumerate}
\item Large enough so that a high compression ratio can be attained with low reconstruction error.
\item Small enough so that the autoencoder runs quickly.
\end{enumerate}
The second item poses a particular challenge in the present case as the DG solver
used to solve the forward acoustic wave equation is highly tuned. Special
care must be taken in order to develop a method that is faster than checkpointing
while maintaining sufficient accuracy. We found empirically that an input dimension of
4096 adequately balanced these two constraints. This corresponds to 8 elements
per MPI rank with 3$^{\text{rd}}$ order polynomials, resulting in 64 nodes per element ($4^3$). 

\subsection{Compression on non-uniform mesh}
\label{section:time_compression}
While using an autoencoder to compress in space is intuitive and simple to implement,
there are severe limitations to the practical applications of such an approach.
First, each MPI rank must have the exact same number of degrees of freedom that
are arranged in exactly the same manner
in order to avoid padding or communicating ghost elements.
This is because autoencoders that employ dense hidden layers have a fixed input size.
Unless extreme care is taken, this requires that the mesh be uniformly refined
and have no hanging nodes.
Such a restriction poses major challenges for problems with a non-rectangular domain
or when non-hexahedral elements are desired.
To mitigate these limitations, we propose to compress along the time dimension, 
treating each element as an independent state.
This method works seamlessly with h-adaptive mesh refinement (refining the mesh geometry), though
p-adaptivity (changing the polynomial order of the element) remains a challenge.
For this work, we only consider algorithms with h-adaptivity.

\begin{figure}[htb!]
  \centering
  \begin{subfigure}{0.5\textwidth}
    \includegraphics[trim=1cm 1cm 4.5cm 1cm,clip=true,width=\textwidth]{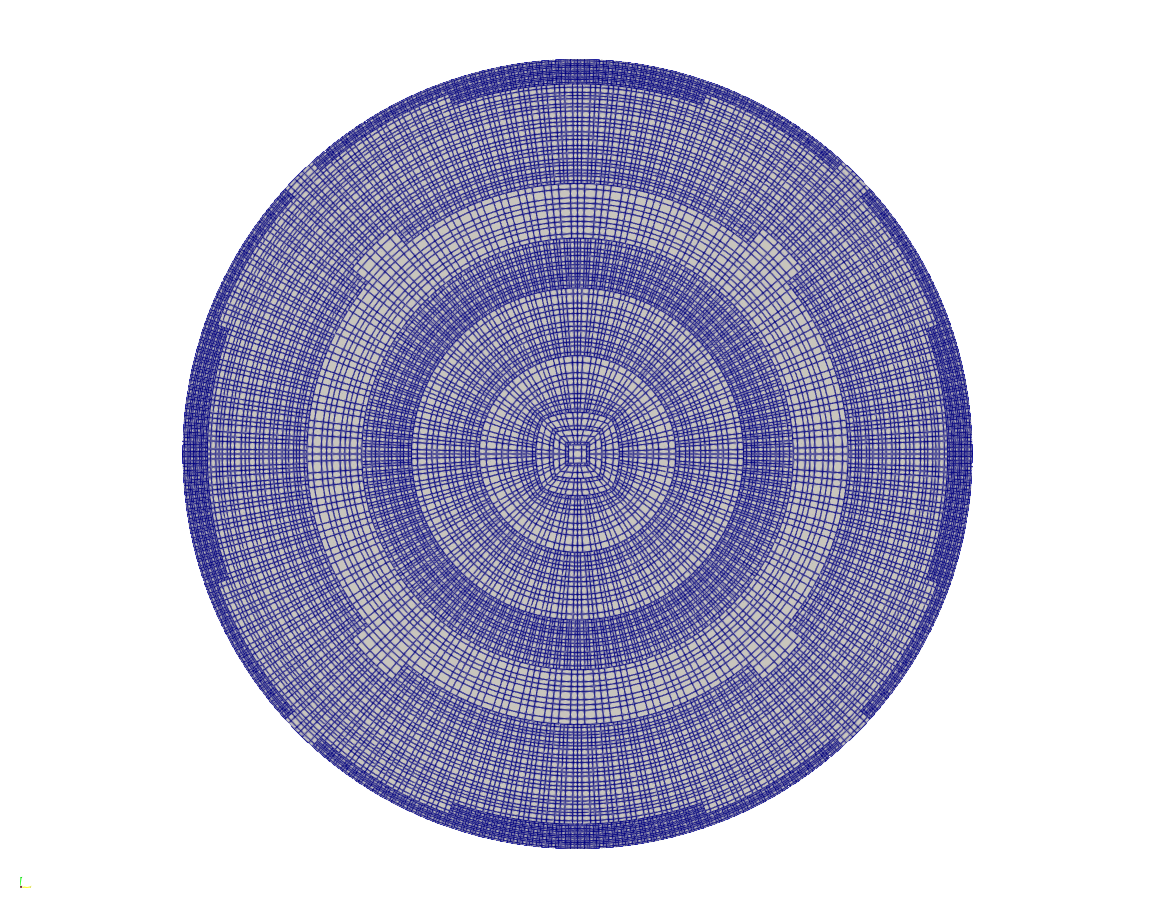}
  \end{subfigure}
  \caption[Earth domain mesh]{
    A cross-section of the earth domain with h-adaptive mesh refinement. Notice the non-uniform
    character of the mesh, with very small elements near the crust required to
    resolve the seismic waves. 
  }
  \label{fig:earth_mesh}
\end{figure}

To demonstrate the success of this approach, we consider solving a seismic inverse problem on
a spherical domain. We employ
h-adaptive mesh refinement to ensure that the seismic waves can be resolved on each element
with at least 3 nodes per wavelength. Figure \ref{fig:earth_mesh} shows a cross section of the
mesh. Hanging nodes that result from local mesh refinement pose no problem for
our proposed compression strategy since each node is considered independently.

The same constraints as discussed in Section \ref{section:space_compression} apply in this
scenario. Similar to the box domain case, we found that an input dimension to the
autoencoder of 4096 worked well. Since we use a 4$^{\text{th}}$ order Runge-Kutta (RK)
explicit timestepping scheme \cite{quarteroni2010numerical}, a 4096 dimensional input
vector results from using $3^{\text{rd}}$ order polynomials (64 nodes per element) and 16 timesteps,
each with 4 RK stages. 

\subsection{Autoencoder architecture and training details}
Any discussion of results obtained using deep learning would be incomplete without
a discussion of the architecture used to obtain the results.
While the order of the data differs between space compression and time compression methods, 
the same autoencoder architecture can be used in either case since the input
and output dimensions are the same. The difference between the trained autoencoders
used for each compression method lies solely in the training data.
Table \ref{table:architecture} shows the architecture design of both the encoder
and the decoder. All layers are ``dense'' layers, which corresponds to an affine
transformation followed by a nonlinear activation.
We found the exponential linear unit
(ELU) to perform better than the rectified linear unit (ReLU), hyperbolic tangent (tanh),
and sigmoid activations. The ELU activation function is given by
\begin{align*}
  \text{ELU} (y) :=
  \begin{cases}
    e^{y} - 1 & \text{if } y < 0 \\
    y & \text{else}.
  \end{cases}
\end{align*}
Table \ref{table:hyperparameters} shows the relevant hyperparameters
of the training process. Training took approximately 3 hours using two Nvidia
2080Ti GPUs. 
\newcommand{\spacing}{6pt}
\begin{table}[htb!] 
  \centering
  \begin{tabular}{l | c | c}
    & Encoder & Decoder \\\hline
    &&\\[-0.2cm]
    Input dimension & 4096 & 64\\[\spacing]
    Output dimension & 64 & 4096\\[\spacing]
    \# hidden layers & 7 & 7 \\[\spacing]
    \multirow{2}{3cm}{Neurons per hidden layer} &
    \multirow{2}{3cm}{\centering [512, 256, 256, 256, 128, 64, 64]} &
    \multirow{2}{4cm}{\centering [128, 128, 256, 256, 256, 512, 4096]}\\&&\\[4pt]
    Activation & ELU & ELU 
  \end{tabular}
  \vspace*{0.2cm}
  \caption[Autoencoder architecture]{
    Autoencoder architecture and hyperparameters used to train network.}
  \label{table:architecture}
\end{table}

\begin{table}[htb!] 
  \centering
  \begin{tabular}{l | c }
    Precision & float32 \\[\spacing]
    Optimizer & LAMB \cite{you2019large} \\[\spacing]
    Initial learning rate & $10^{-3}$ \\[\spacing]
    \multirow{2}{*}{Learning rate decay} & $lr = 0.5 lr$ \\[-2pt]
              & applied every 5 epochs \\[\spacing]
    Epochs & 20 \\[\spacing]
    Batch Size & 512 \\[\spacing]
    \# GPUs &   2$\times$ 2080Ti
  \end{tabular}
  \vspace*{0.2cm}
  \caption[Training hyperparameters]{
    Training hyperparameters used to fit autoencoder weights and biases to data. }
  \label{table:hyperparameters}
\end{table}

As mentioned in Section \ref{section:space_compression}, it is difficult to
design an autoencoder compression system that both performs well in terms of accuracy and compression ratio
while also being fast.
Recall that the evaluation of a single dense layer requires the computation of
\begin{equation}
  \label{eqn:single_layer}
  \sigma \LRp{\bs{y} \bs{W} + \bs{b}}. 
\end{equation}
Considering the first layer of the encoder, $\bs{W}_1 \in \R^{4096 \times 512}$, $\bs{y} \in \R^{4096}$
and $\bs{b}_1 \in R^{512}$. There are 2,097,152 entries in the matrix $\bs{W}_1$. Compare this
to the second layer which has a weight matrix $\bs{W}_2 \in \R^{512 \times 256}$ and 131,072 entries.
The first layer then has 16 times more entries in the weight matrix than the second layer.
Further layers of the encoder have even smaller weight matrices. It is clear that the
dominant cost in compressing a vector is the computation of the first matrix-vector multiply,
$\bs{y} \bs{W}_1$.
To mitigate this cost, we develop a hybrid sparse-dense architecture whereby the first layer 
of the encoder and last layer of the decoder are trained to be 95\% sparse.

As an example to show the improvement of the sparse-dense architecture over the
fully dense architecture, consider compressing a \textit{batch} of 864 vectors
each with size 4096. Since evaluating the output of the neural network requires nothing
more than matrix-vector multiplications, vector additions, and element-wise
application of an activation function, all 864 vectors can be processed
simultaneously by stacking into rows.
Then the input $\bs{y}$ has shape $(864, 4096)$.
Timing results for both the encoder and decoder
on a single core of an Intel Xeon Platinum 8280 CPU on Frontera
are shown in Table \ref{table:autoencoder_timings}. We would then expect approximately
a 50\% speedup with the sparse-dense architecture over the fully dense architecture.
\begin{table}[htb!] 
  \centering
  \begin{tabular}{l | c | c}
    & Encoder & Decoder \\\hline&&\\[-8pt]
    Fully dense &         102 ms  & 140 ms \\[\spacing]
    Hybrid sparse-dense & 40 ms & 79 ms
  \end{tabular}
  \vspace*{0.2cm}
  \caption[Single core timing]{
    Comparison of timings of fully dense architectures vs hybrid sparse-dense architecture
    for encoder and decoder for compressing/decompressing 864 vectors
    on a single core of an Intel Xeon Platinum 8280 CPU.  }
  \label{table:autoencoder_timings}
\end{table}

Unfortunately, the application setting is slightly different since
all cores are used rather than just a single core. Consider the same test
case replicated on each CPU core simultaneously. There are 56 cores
across 2 sockets on a single node of Frontera. 
Since each task executes independently with no communication, 
we might expect perfect scaling and for the test to complete in the same amount of time.
However, there are now more cores contending for data with limited cache and memory bandwidth.
Indeed, in the first case, 864 single precision input vectors require only 14.2 MB
to store. This easily fits in the 38.5 MB of L3 cache of a single Intel Xeon Platinum 8280 \cite{intel}.
Furthermore, the weights of the fully dense encoder require only 17.8 MB to store.
The input data and the entire encoder can fit in L3 cache! Since L3 cache is shared among
all cores of this CPU, bytes must be streamed from main memory when all cores are
running simultaneously.
Table \ref{table:autoencoder_timings_56} shows the
timings for the encoder and decoder with the fully dense and hybrid sparse-dense architectures
when running on all cores.
\begin{table}[htb!] 
  \centering
  \begin{tabular}{l | c | c}
    & Encoder & Decoder \\\hline&&\\[-8pt]
    Fully dense &        110  ms &  160 ms \\[\spacing]
    Hybrid sparse-dense & 120 ms &  117 ms
  \end{tabular}
  \vspace*{0.2cm}
  \caption[Multi-core timing]{
    Comparison of timings of fully dense architectures vs hybrid sparse-dense architecture
    for encoder and decoder for compressing/decompressing 864 vectors
    on 56 cores of 2 Intel Xeon Platinum 8280 CPUs on a single node of the Frontera supercomputer.  }
  \label{table:autoencoder_timings_56}
\end{table}

Although there is still significant advantage in using the hybrid architecture for the decoder,
the evaluation of the encoder proves to be faster using the fully dense architecture than the
hybrid architecture when all 56 cores are being utilized.
We are unsure why this is the case at the time of writing
and further optimization is a topic for future investigation.
However, measuring the performance showed that
we could achieve faster compression by using a fully dense matrix for the weights of the first layer
while optimizing the evaluation time of the decoder using a sparse matrix for the final layer. Because of the strict speed 
requirements, we do not utilize any error bounding methods
which require the computation and potential storage 
of a residual vector \cite{lee2022error,liu2021exploring,liu2021high,abu2014efficient}.

\subsubsection{Training data}
\label{section:training_data}
One common criticism of deep learning methods is that 
large training datasets are often required
for the DNN to perform well in practice.
In the case of compression, other algorithms such
as those implemented in the ZFP package \cite{lindstrom2014fixed}
do not require training data. 
Even when deep learning methods may perform well on training data,
they may perform poorly in practice when the training dataset is not sufficiently
large, a concept referred to as \textit{generalization error}.
While this is true, we argue that the cost of training an autoencoder is a one-time
upfront cost.
Further, data generation and training require several orders of magnitude
fewer compute node-hours than solving the inverse problem.
In this section, we detail the data generation process.

Consider a statistical interpretation where the training data $\bs{y}\in \R^{n} $, is drawn from
some distribution, $\probMeas[\bs{y}]$.
Suppose that we have a dataset of $N$ samples $\mathcal{D}^N=\LRc{\yb^1,..., \yb^N}$ where
$\yb^i \sim \probMeas[\yb]$.
Let $\Psi_{\mathcal{D}^N}\LRp{\yb}$ denote the output of
the machine learning algorithm after being trained with dataset $\mathcal{D}^N$.
A common metric for evaluating the performance of any machine learning algorithm is the 
\textit{generalization error}.
Given some loss function denoted $\loss\LRp{\Psi\LRp{\yb}, \yb}$,
such as mean squared error,
\begin{equation*}
  \loss_{\text{MSE}}\LRp{\Psi\LRp{\bs{y}_{\text{true}}}, \bs{y}_{\text{true}}}
  := \frac{1}{n}\norm{\bs{y}_{\text{true}} - \Psi\LRp{\bs{y}_{\text{true}}}}^2_2,
\end{equation*}
the generalization error for an autoencoder is
$\expect[\probMeas[\yb]][\loss\LRp{\Psi\LRp{\yb}, \yb}]$
\cite{nadeau1999inference}. That is, the generalization error is the average
loss over the entire distribution $\probMeas[\yb]$.

Generalization error plays an important role in designing a data generation process.
There are two ways one can reduce the generalization error:
\begin{enumerate}
\item Increase the size of the training dataset by drawing more samples from $\probMeas[\yb]$.
\item Choose $\probMeas[\yb]$ so that it can be well-approximated with few samples. 
\end{enumerate}
The second method can be interpreted as choosing the most narrow distribution
that is appropriate for the application at hand.
We will explore how this option can be used in the context of generating data for seismic inversion. For typical machine learning applications, this may not be possible as we don't know $\probMeas[\yb]$. The Bayesian inversion setting makes this possible for our application, as we discuss next.

In order to compress seismic data, we need to generate examples
similar to the states that will be seen in application. 
One way to do this is to solve the acoustic
wave equation and save some of the data for training.
However, it is still unclear which solutions of the wave equation should be included
in the training dataset.
Should we try to train the
autoencoder to compress and decompress all possible solutions to the wave equation on a given domain?
There are clearly wave equation solutions that are not
relevant for our application.
Mathematically, the challenge becomes defining the relevant distribution $\probMeas[\yb]$.
The Bayesian framework \cite{kaipio2006statistical}
helps us to choose a good distribution for generating data. 
Bayes' formula tells us,
\begin{equation}
    \probMeas[post][\qoi \given \obs] \propto \probMeas[like][\obs \given \qoi] \probMeas[prior][\qoi].
\end{equation}
This can be interpreted as the likelihood, $\probMeas[like][\obs \given \qoi]$,
updating the prior density, $\probMeas[prior][\qoi]$, based on observed data, $\obs$,
to produce the posterior density, $\probMeas[post][\qoi \given \obs]$,
\cite{scales2001prior}.
The observations $\obs$ are related to the states
that will be compressed, denoted $\yb$ here, through an observation 
operator, i.e.  $\obs = \bs{\mathcal{B}}\yb$. 
Since the posterior density 
is the product of the likelihood and the prior,
the posterior density is
non-zero only where the prior density is non-zero.
That is, the posterior lies within the support of the prior. 
Therefore, we are only interested in whether the autoencoder
performs well in compressing states, $\yb$, that can be generated by
samples from the prior.
We can generate relevant examples of possible states by solving the acoustic wave
equation with the PoI (the wave speed) set to samples from the prior. 
The autoencoder does not need to perform well on every possible
state that could be generated from the wave equation, it only needs to work well on states
generated by solving the wave equation with parameters sampled from the prior.

With motivation from statistical learning theory and Bayesian inverse problems, 
we now detail an algorithm for generating training data.
\begin{algorithm}[h!t!b!] 
  \caption{Data generation algorithm for training compression autoencoder}
  \hspace*{\algorithmicindent} \textbf{Input}: number of samples, $n$; final time, $T$\\
  \begin{algorithmic}[1]
    \FOR{$\, i=1,...,n$}
    \STATE {$\cp = $ random draw from prior}\\
    \FOR{$\, t = 1, ..., T$}
    \STATE {Solve acoustic wave equation for one timestep}\\
    \STATE {Consolidate states, $\velocity$, $\strain$, into compression data structure}
    \IF {Data structure is full}
    \STATE {Write data structure to file}
    \ENDIF
    \ENDFOR
    \ENDFOR
  \end{algorithmic}
  \label{algorithm:training_data}
\end{algorithm}
The compression data structure is a vector that stores the state variables
in the order corresponding to the compression scheme in use.
This setup allows any ordering to be used  and for training data
to automatically have the correct ordering --- reducing effort to set up the training pipeline. 
Algorithm \ref{algorithm:training_data} can generate large amounts of data very quickly.
For the earth-scale problem with an 11.2 million DoF mesh,
several terabytes can be generated in only a few minutes.
To reduce the amount of data written to file, we randomly chose whether to write or discard the data structure.
On 32 nodes of Frontera, solving the wave equation and writing to file takes on the order of 1 minute.
With 10 draws from the prior, keeping only $0.1\%$ of the data, we generated 50 GB of training
data in less than 1 hour.
The total cost of data generation is then $\bigO \LRp{10}$ node-hours. 

Finally, we discuss normalization of the training data.
The states have scales that vary between $10^{-16}$ and $10^{5}$ over the
solution of \eqref{eqn:model}.
It would be difficult for an autoencoder to
be able to accurately compress and decompress states with such widely varying scales,
especially in single precision.
To reduce the burden on the autoencoder of learning scale, we normalize the data between 0 and 1.
This enables us to store the scale and offset of the data in double precision while capturing the
relative variation of the data in single precision.
Consider an arbitrary state vector that will be compressed,
denoted $\yb$ following the notation in \eqref{eqn:single_layer}. Let $\yb_{\text{true}}$ be
a consolidated state vector that will be compressed
(state DoFs consolidated either in space or in time).
The input vector to the encoder, $\yb$, is then given by 
\begin{align*}
  r   &= \max(\max(\yb_{\text{true}}) - \min(\yb_{\text{true}}), \beta) \\
  \yb &= \LRp{\yb_{\text{true}} - \min(\yb_{\text{true}})} / r.
\end{align*}
In order to avoid division by zero, $r$ is set to be at least some tolerance, $\beta$.
We choose $\beta = 10^{-7}$.

\section{Numerical results}
\label{section:ae_numerical_results}
In this section we present numerical results showing the efficacy of our proposed
autoencoder compression approach for both spatial compression on the box domain problem
(Section \ref{section:space_compression}) and temporal compression on a spherical domain
(Section \ref{section:time_compression}).
For both problems, we solve the inverse problem using a  Newton conjugate gradient
(NewtonCG) method \cite{yang2009newton,epanomeritakis2008newton}. In each conjugate gradient iteration, the action of the Hessian
on a vector is required which entails solving two additional linearized
PDEs, referred to as the incremental
forward problem and incremental adjoint problem.
Expressions for these PDEs in the case of the acoustic wave
equation can be found in \cite{bui2012extreme}. 
The incremental forward problem relies on the solution of the forward problem
at each timestep and so requires decompression as well. The incremental adjoint problem
requires the solution of both the forward and incremental forward problems,
and thus requires another decompression of the forward states.
Empirically, we find that the gains of decompression are mostly offset by the
costs of compression. However, the states only need to be compressed once
per Newton iteration and are potentially decompressed many times.
That is, the cost of compression is amortized over every evaluation of the adjoint,
incremental forward, and incremental adjoint problems within a single Newton iteration.
While it is also possible to compress and decompress the incremental forward problem,
it is not explored in this paper since each incremental forward solution is only re-used
once. 

\subsection{Results: compression in space}
We demonstrate the compression in space approach on a cube domain with
all sides of length 1km. Five sources are located at 1m depth 
and 1089 receivers are evenly spaced at the surface,
observing the velocity field for 3 seconds.
The sources are Ricker wavelets in time, smoothed by convolving
with a narrow Gaussian in space.
The source term for a source at $\xb_0$ is then given by 
\begin{equation*}
  \bs{g}\LRp{\xb, t; \xb_0} = \rho\LRp{\xb} \frac{1}{\sqrt{2\pi}\sigma_{\xb}}
  e^{- \frac{\norm{\xb - \xb_0}_2^2}{2 \sigma_{\xb}^2}}
  \LRp{1 - \frac{(t - t_c)^2}{\sigma_t^2}}e^{-\frac{(t-t_c)^2}{2\sigma_t^2}}
  \LRs{0, 0, -1}^T.
\end{equation*}
In this case, we consider a wavelet
centered at $t_c = 0.6$ seconds with $\sigma_t = \frac{1}{\pi}$
and $\sigma_{\xb} = 0.05$.

\begin{figure}[htb!]
  \centering
  \begin{subfigure}{0.45\textwidth}
    \includegraphics[trim=12cm 12cm 12cm 10cm,clip=true,width=0.95\textwidth]{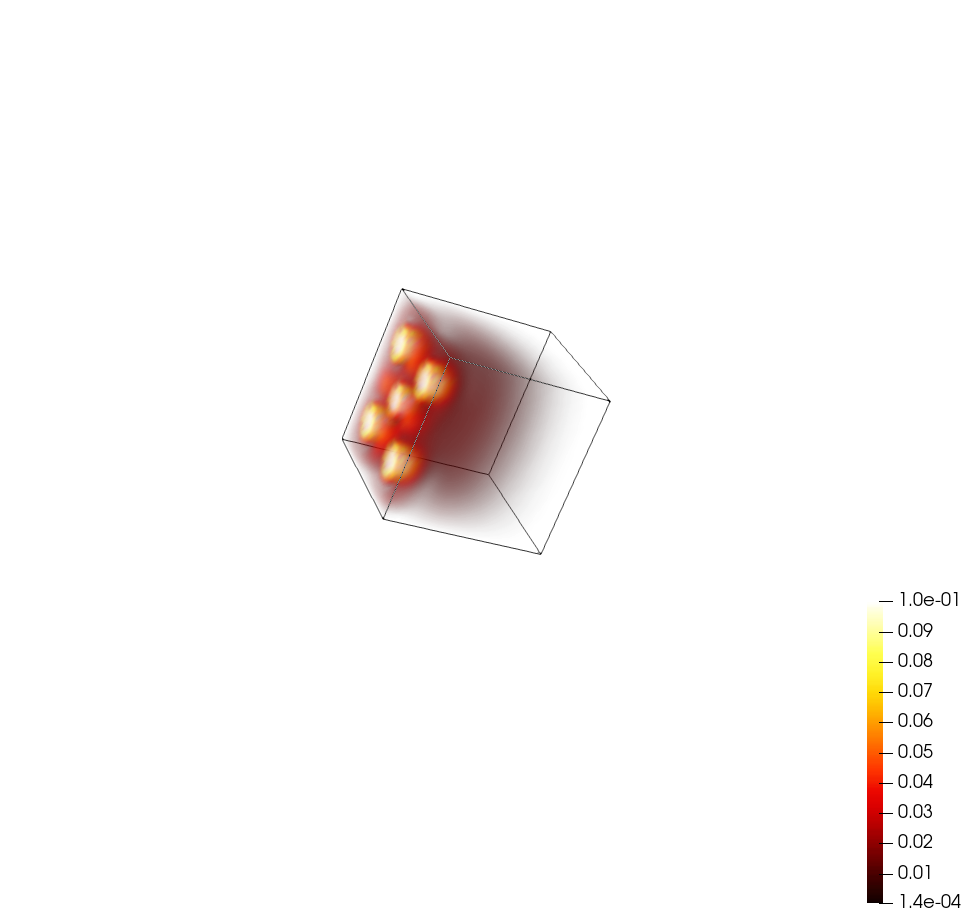}
    \caption{}
    \label{fig:box_state}
  \end{subfigure}
  \begin{subfigure}{0.45\textwidth}
    \includegraphics[trim=8cm 2.5cm 10cm 2cm,clip=true,width=\textwidth]{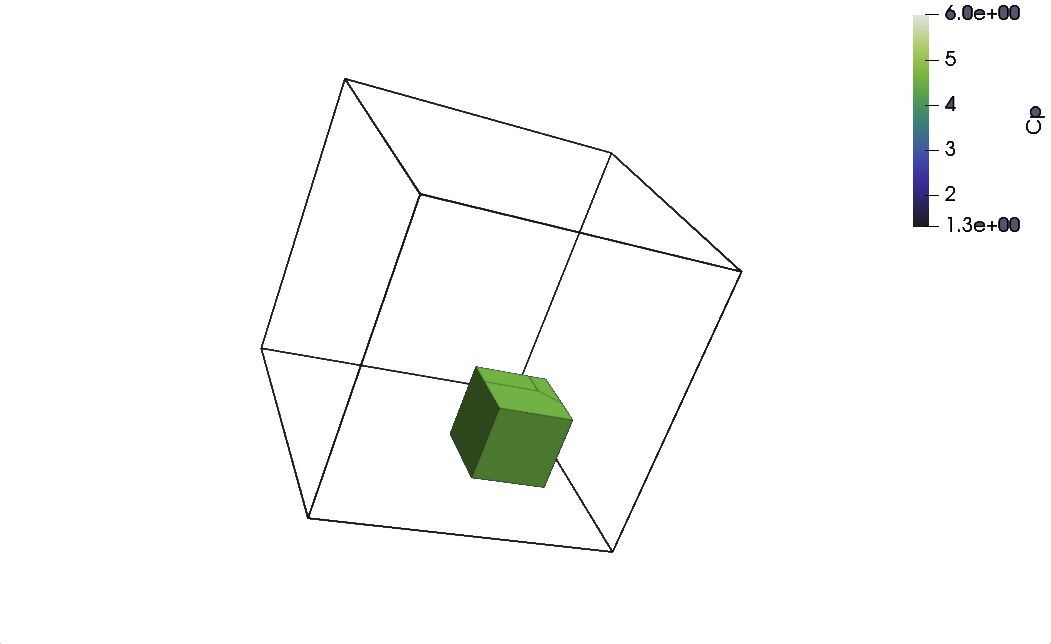}
    \caption{}
    \label{fig:box_true_solution}
  \end{subfigure}
  \caption[Box state and true parameter]{
    (a) Snapshot of velocity magnitude at 1.2 seconds
    (b) Inclusion in background wave speed.}
\end{figure}

A snapshot of the velocity magnitude is shown in Figure \ref{fig:box_state}.
Following the setup in Section \ref{section:space_compression}, we show
both speedup and absolute error results compared to the checkpointing solution.
While the synthetic true solution is our target, our compression
strategy aims to replace checkpointing as a means for managing limited memory, 
not to change the inverse solution.
Therefore, the
appropriate comparison is with the inverse solution 
computed using checkpoints rather than the true solution.
Similar to the earth problem, we invert for a variation in the background
wavespeed field.
A box-shaped inclusion is introduced into the
domain, shown in Figure \ref{fig:box_true_solution}, and is referred to as the anomaly.

\begin{figure}
  \centering
  \begin{subfigure}{0.43\textwidth}
    \includegraphics[trim=10cm 3cm 10cm 6cm,clip=true,width=\textwidth]{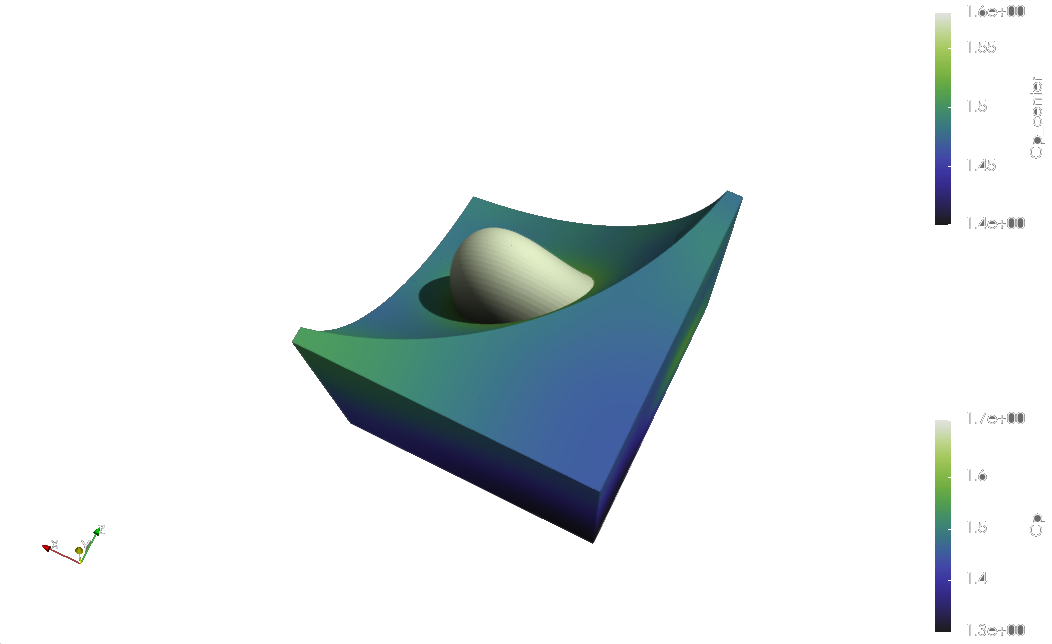}
    \caption{Checkpoints}
  \end{subfigure}%
  \begin{subfigure}{0.43\textwidth}
    \includegraphics[trim=10cm 3cm 10cm 6cm,clip=true,width=\textwidth]{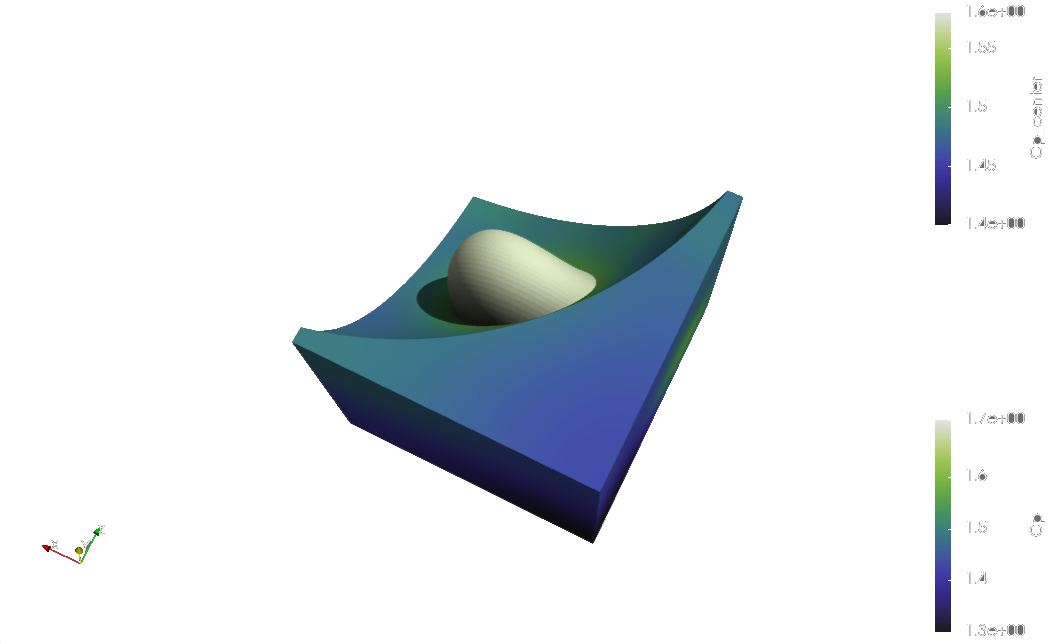}
    \caption{Autoencoder compression}
  \end{subfigure}%
  \begin{subfigure}{0.15\textwidth}
    \includegraphics[trim=1.1cm 0cm 1.8cm 0.6cm,clip=true,width=0.55\textwidth]{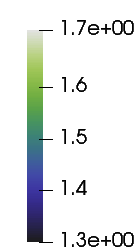}
    \caption{$\cp$, km/s}
  \end{subfigure}%
  \caption[Box inverse solution]{
    Inverse solutions on box domain using (a) checkpoints and
    (b) autoencoder compression in space for problem with 16,777,216 DoF mesh.
    The solutions have the same color scale with the inclusion rendered as
    an isosurface at 1.67 km/s.
  }
  \label{fig:uqbox_cutaway}
\end{figure}

To demonstrate the weak scaling of this method, we show results for several
refinement levels, each refinement increasing the total number of degrees of freedom by 8 times.
Table \ref{table:box_speedup} shows the speedup compared to the checkpointing solution.
While the overall speedup may not be very impressive, it is not the only objective of our proposed compression approach. As we have discussed, the other objectives are to reduce memory footprint and to eradicate checkpointing. 
Furthermore, there is a maximum speedup of 1.3$\times$
for the gradient since the time to solve
the forward problem, the adjoint problem, and evaluate the gradient are approximately equal.
That is, the cost of computing the gradient in the ideal case is 3 PDE solves while the checkpointing
case requires 4 PDE solves. Thus, the speedup of 1.22 for the largest mesh in Table \ref{table:box_speedup} is close to ideal. Similarly, the Hessian-vector product computation has the maximum
speedup of 1.5 (requiring 4 PDE solves in the ideal case vs 6 PDE solves with checkpointing), and the speedup of 1.27, without compressing the incremental forward solve, for our compression approach is already an achievement. 
Furthermore, the acoustic wave equation, \eqref{eqn:model},
is linear and not difficult to solve with explicit time stepping.
Our proposed approach would show closer to ideal speedup when used to
replace implicit solvers or when used on nonlinear problems with much more time consuming
forward solves since the cost of autoencoder compression only depends on the number of parameters
of the neural network.
Table \ref{table:box_speedup} shows that the speedup approaches the ideal one and the relative error decreases as the number of DoFs increases.
A cutaway of the inverse solution with the largest mesh for both the checkpointing case and the
autoencoder compression case is shown in Figure \ref{fig:uqbox_cutaway}. 
The fact that the anomaly reconstructions are visibly indistinguishable is due to 0.4\% difference.

\begin{table}[htb!]
    \centering
    \begin{tabular}{l|c|c|c|c}
      Mesh DoFs & Cores &$\grad$ Speedup & $\mathcal{H}\bs{v}$ Speedup & Relative $l_2$ Error \%  \\[2pt]
      \hline &&&  \\[-8pt]
      4,096      & 1    & 1.06   & 1.20   & 1.7   \\[6pt] 
      32,768     & 8    & 1.19   & 1.26   & 0.7   \\[6pt] 
      262,144    & 64   & 1.17   & 1.22   & 0.6   \\[6pt] 
      2,097,152  & 512  & 1.17   & 1.23   & 0.7   \\[6pt] 
      16,777,216 & 4096 & 1.22   & 1.27   & 0.4   \\[6pt] 
    \end{tabular}
    \caption[Box domain results]{
      Speedup of the autoencoder compression algorithm compared to the checkpointing solution.
      $\grad$ represents the gradient and $\mathcal{H}\bs{v}$ represents a Hessian-vector product. 
      The speedup approaches the ideal speedup as the number of DoFs increases while the relative error decreases. }
    \label{table:box_speedup}
\end{table}

\subsection{Results: compression in time}
\begin{figure}[htb!]
  \centering
  \includegraphics[width=0.45\textwidth]{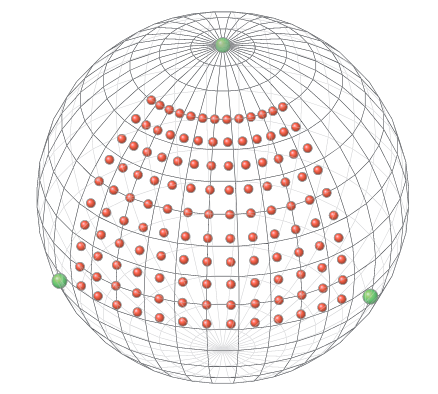}
  \caption[Earth sources and receivers]{
    Three sources shown in green at the North Pole and along the equator.
    130 receivers shown in red distributed in North America
    spaced 7.5 degrees apart.
    Figure from \cite{bui2013computational}, 
    Copyright \copyright2013 Society for Industrial and Applied Mathematics. Reprinted with permission. All rights reserved.
  }
  \label{fig:earth_source_receiver}
\end{figure}

While it is intuitive that seismic waves exhibit strong spatial redundancy
and thus have high compressibility, it is not immediately obvious that
temporal compression schemes will yield acceptable results.
We test the temporal compression scheme presented in Section \ref{section:time_compression}
on the earth-scale seismic inverse problem.
We build upon the framework established in \cite{bui2012extreme,bui2013computational}
to show that our proposed compression scheme is a viable method for
large-scale seismic inverse problems. The results shown here
follow the setup of Case \textit{II} from \cite{bui2013computational}.
There are 3 source locations in the Northern Hemisphere,
one at the North Pole and two placed at 90 degree increments along the equator.
Here, we consider sources modeled as Dirac delta functions in space centered
at the source location, $\xb_0$, 
and a Gaussian in time given by
\begin{equation*}
  \bs{g}\LRp{\xb, t; \xb_0} = \frac{1}{\sqrt{2\pi}\sigma}
  \delta \LRp{\xb - \xb_0} e^{-\frac{\LRp{t - t_c}^2}{2\sigma^2}}
\end{equation*}
where $t_c = 60$ seconds and $\sigma = 20$ seconds. 
We place 130 receivers spaced at 7.5 degree increments throughout North America
as show in Figure \ref{fig:earth_source_receiver}, measuring
the velocity field.
Additive noise with standard deviation 0.002 is introduced into the data.
The earth is modeled as a sphere with radius 6,371km.
Both the sources and receivers are buried at 10km depth. 

Course meshes cannot resolve the sufficiently high frequency waves
required to effectively image the earth and thus may not be used. 
We found that a maximum frequency resolution of 0.05Hz
was sufficient to obtain good inversion results.
Thus, the smallest domain size that we consider has a mesh with
11.2 million degrees of freedom. This corresponds to a discontinuous Galerkin
discretization with 174,481 elements using 3$^{\text{rd}}$ order polynomials.
We solve this problem on 64 nodes (3,584 CPU cores) of Frontera
running for 10 hours. Speedup results compared to the checkpointing case
are shown in Table \ref{table:uqearth_speedup}.
Here we report the average relative $l_2$ error of the solution obtained using autoencoder compression
to the solution obtained using checkpointing. 
Similar to the box domain, there is greater speedup in the $\mathcal{H}\velocity$ 
product computation than the gradient computation due to paying the cost 
of compression during the gradient solve. 
The difference in speedup between the box problem and the earth problem 
is likely due to variable per-process problem size induced by the non-uniform 
spherical mesh and related cache contention issues. 
As mentioned in Section \ref{section:ae_methodology}, 
cache effects can greatly impact the performance of the compression algorithm.

\begin{table}[htb!]
    \centering
    \begin{tabular}{l|c|c|c}
      Mesh DoFs & $\grad$ Speedup & $\mathcal{H}\bs{v}$ Speedup & Relative $l_2$ error \%  \\[2pt]
      \hline &&&  \\[-8pt]
      11,166,784 &  1.05  & 1.17  &  0.03
    \end{tabular}
    \vspace*{0.2cm}
    \caption[Earth domain results]{
      Speedup of the autoencoder compression algorithm compared to the checkpointing solution.
      $\grad$ represents the gradient and $\mathcal{H}\bs{v}$ represents a Hessian-vector product. 
      Timing results were obtained by running on 64 nodes (3,584 CPU cores) of the Frontera supercomputer. 
    }
    \label{table:uqearth_speedup}
\end{table}

To demonstrate scalability of the autoencoder compression approach,
we solve the same problem
on a refined mesh with 139,227,136 DoFs. The full state vector (3 velocities
and 3 strains) then has 835,362,816 DoFs.
This problem was solved on 256 nodes of Frontera
(14,336 CPU cores).
Table \ref{table:zfp_vs_ae_large} compares the autoencoder compression
results to the checkpoint solution. 
We find that there is a slight reduction in terms of gradient speedup and a
slight increase in speedup of the Hessian-vector product compared to
the course mesh solution detailed in Table \ref{table:zfp_vs_ae}.
One possible reason for the decrease in speedup of
the gradient computation is the increased per-process work of the larger problem. 
The smaller problem with 11 million mesh DoFs has an average of 
4,681 mesh DoFs per process, while the larger problem has 
an average of 9,711 mesh DoFs per process. 
We explored the effects of limited cache size on the timing 
of autoencoder compression in Section \ref{section:ae_methodology} 
and found that the compression is more sensitive to 
cache effects than decompression, at least for our implementation. 
It is therefore not surprising that the speedup of the gradient 
computation, which includes the cost of compression, is reduced 
when each process is responsible for compressing a larger amount of data. 
While speedup is desired in all computations,
there is an order of magnitude more Hessian-vector products
than gradient computations required to compute the MAP solution,
with $\bigO\LRp{100} \mathcal{H}\bs{v}$ products and
$\bigO\LRp{10}$ gradient evaluations. The overall cost is dominated
by $\mathcal{H}\bs{v}$ products, and so the observation that the
gradient computation sees no speedup from compression
has little impact on the total computation time. 

\begin{table}[htb!]
  \centering
  \begin{tabular}{p{3cm} | p{2cm}|p{2.2cm}|p{2.3cm}}
    Mesh DoFs& $\grad$ Speedup & $\mathcal{H}\bs{v}$ Speedup &
      Relative error \% \\[2pt] \hline
    &&& \\[-8pt]
    139,227,136 & 1.00  & 1.18  &  0.06\\[6pt]
  \end{tabular}
  \vspace*{0.2cm}
  \caption[Earth domain results - large]{
    Inverse results for earth-scale seismic inverse problem with higher mesh refinement,
    solved on 14,336 CPU cores. 
    The gains from decompression over resolving from checkpoints for computation of the
    gradient is completely offset by the cost of compression. However, this
    cost is paid once and there is similar speedup of the Hessian-vector product
    as seen on a smaller problem in Table \ref{table:zfp_vs_ae}
  }
  \label{table:zfp_vs_ae_large}
\end{table}

\begin{figure}[htb!]
  \centering
  \begin{subfigure}{0.3\textwidth}
    \includegraphics[trim=4.7cm 4.7cm 4.7cm 4.7cm,clip=true,width=\textwidth]{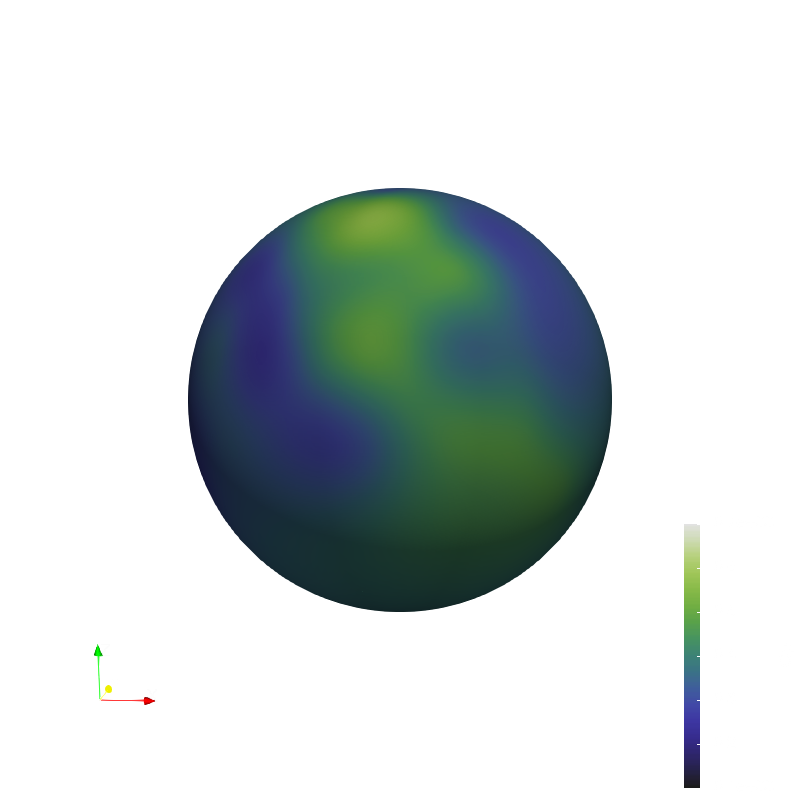}
    \caption{Checkpoint}
  \end{subfigure}%
  \begin{subfigure}{0.3\textwidth}
    \includegraphics[trim=4.7cm 4.7cm 4.7cm 4.7cm,clip=true,width=\textwidth]{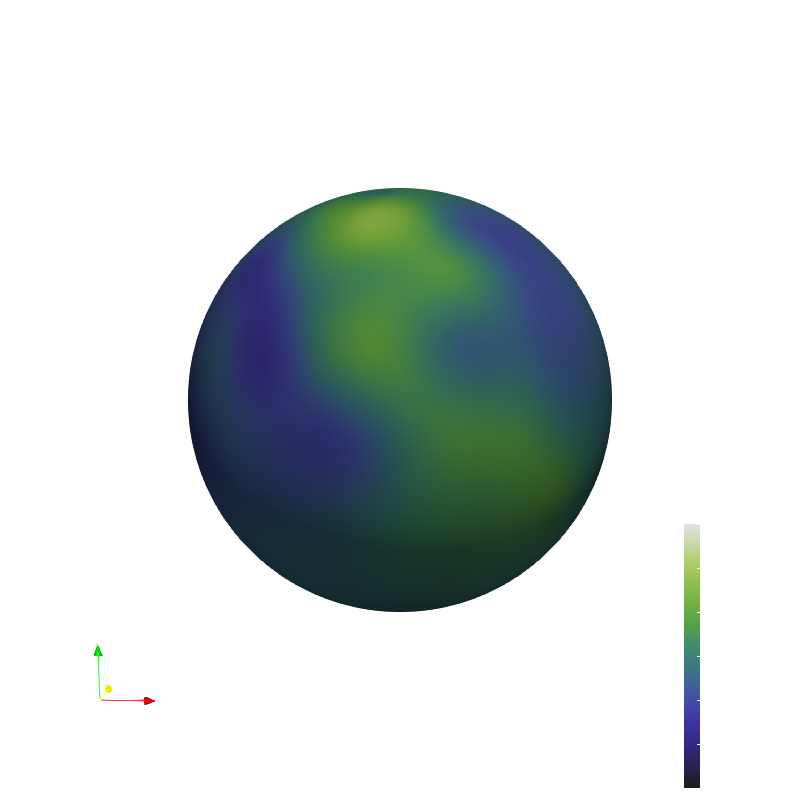}
    \caption{Autoencoder}
  \end{subfigure}%
  \begin{subfigure}{0.3\textwidth}
    \includegraphics[trim=4.7cm 4.7cm 4.7cm 4.7cm,clip=true,width=\textwidth]{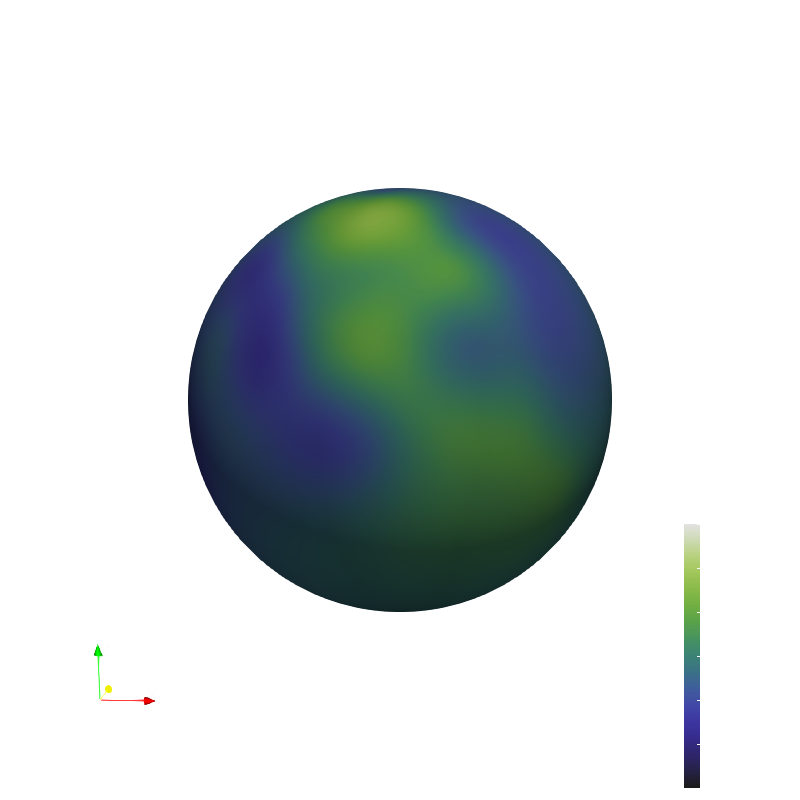}
    \caption{zfp, $\eta = 10^{-5}$}
  \end{subfigure}\\
  
  \begin{subfigure}{0.3\textwidth}
    \includegraphics[trim=4.7cm 4.7cm 4.7cm 4.7cm,clip=true,width=\textwidth]{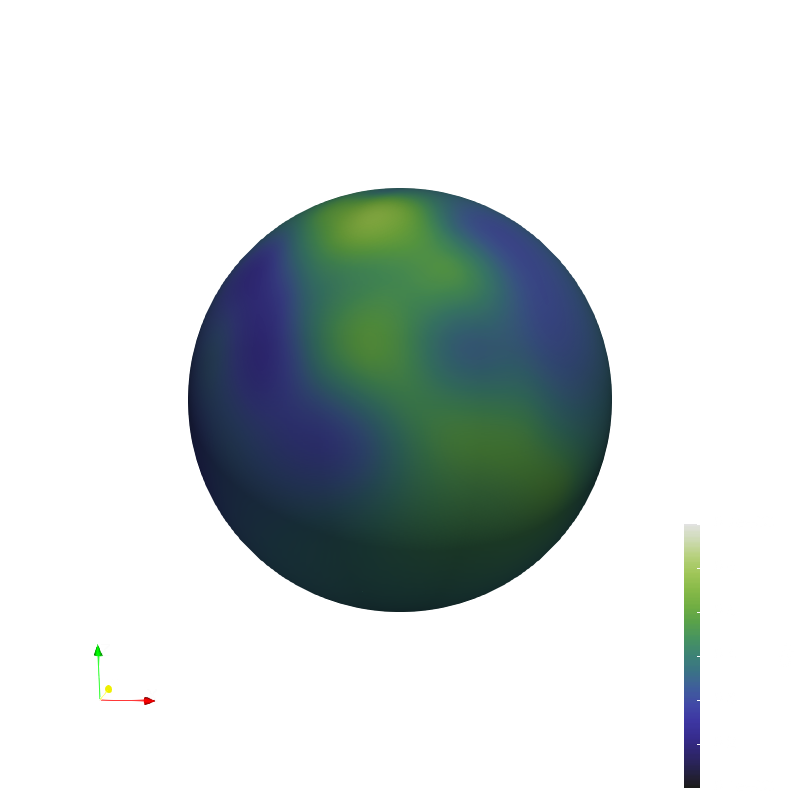}
    \caption{zfp, $\eta = 10^{-4}$}
  \end{subfigure}%
  \begin{subfigure}{0.3\textwidth}
    \includegraphics[trim=4.7cm 4.7cm 4.7cm 4.7cm,clip=true,width=\textwidth]{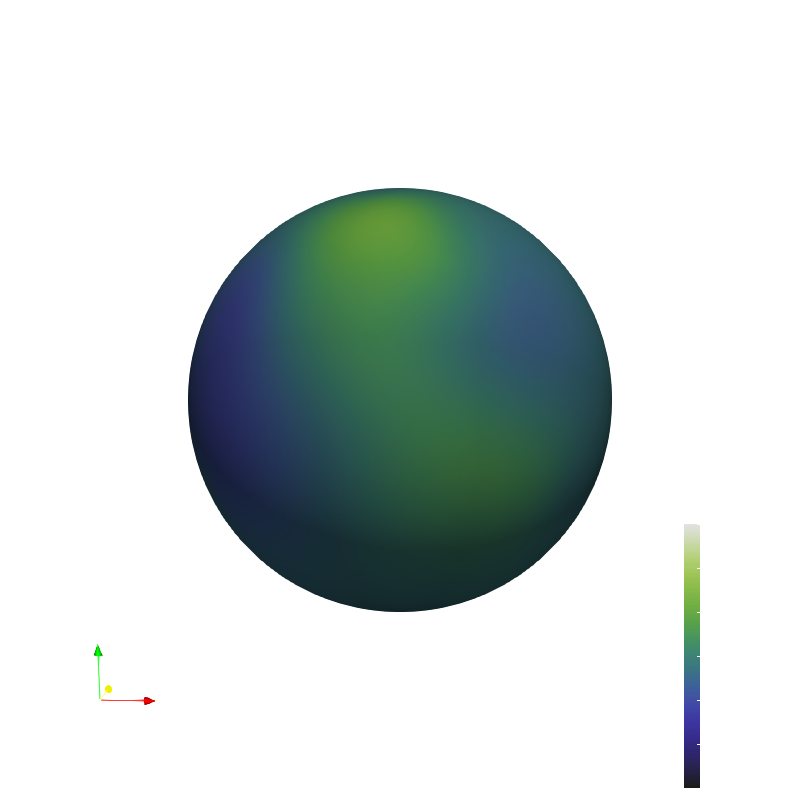}
    \caption{zfp, $\eta = 10^{-3}$}
  \end{subfigure}%
  \begin{subfigure}{0.3\textwidth}
    \includegraphics[trim=4.7cm 4.7cm 4.7cm 4.7cm,clip=true,width=\textwidth]{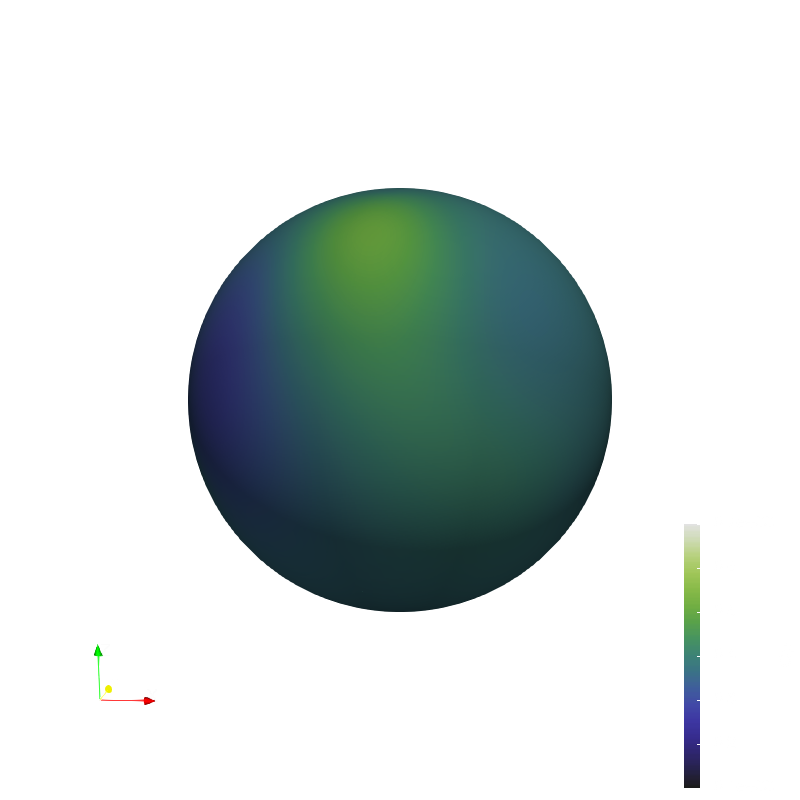}
    \caption{zfp, $\eta =  10^{-2}$}
  \end{subfigure}
  \caption[Comparison to ZFP]{
    Comparison of inverse solutions obtained using autoencoder compression and the ZFP
    compression package with various compression tolerances, $\eta$. Results are shown
    at 130km depth on the same color scale. 
  }
  \label{fig:zfp_vs_ae}
\end{figure}

Having shown that autoencoders are a viable compression technique for accelerating seismic inverse problems,
we now compare our results to the state-of-the-art off-the-shelf
compression technique proposed for seismic inversion in \cite{kukreja2019combining}.
While other off-the-shelf packages have been proposed, the most comprehensive results in the literature are given for the ZFP compression package \cite{lindstrom2014fixed},
especially in the seismic inversion literature \cite{kukreja2019combining,kukreja2022lossy}.
Therefore, we compare the speed, accuracy, and
compression ratio of our autoencoder approach to the ZFP compression package.
While the ZFP package does not require normalization, we tried both with and without normalization
and found that normalizing the states before compressing was necessary to obtain accurate results.
Table \ref{table:zfp_vs_ae} shows the inversion results using ZFP compression for several
different required accuracies. The relative errors for all methods are still quite small,
but ZFP results have qualitatively different solutions than the checkpointing solution
when the compression tolerance $\eta = 10^{-2}$ and $\eta = 10^{-3}$
as shown in Figure \ref{fig:zfp_vs_ae}.
Thus, $\eta \leq 10^{-4}$ is required in order to obtain comparable quality compression
results to the autoencoder solution. 
While autoencoder compression does not show a significant
advantage over ZFP in terms of speed, there is a substantial difference in the
compression ratio, and hence memory footprint. With this tolerance, the compression ratio
of ZFP is 8.4.
On the other hand, the autoencoder compression algorithm
achieved a compression ratio of $128$,
15 times higher than ZFP while maintaining comparable accuracy.

\begin{table}[htb!]
    \centering
    \begin{tabular}{p{3cm}|p{1.8cm}|p{1.8cm}|p{1.6cm}|p{2cm}}
      Method & $\grad$ Speedup & $\mathcal{H}\bs{v}$ Speedup &
      Relative $l_2$ error \% & Compression Ratio \\[2pt] \hline
       &&&  \\[-8pt]
      Autoencoder          &  1.05  & 1.17  &  0.03  & 128 \\[6pt]
      ZFP, $\eta = 10^{-5}$ &  1.0   & 1.09  &  0.004 & 7.2 \\[6pt]
      ZFP, $\eta = 10^{-4}$ &  1.05  & 1.16  &  0.02  & 8.4 \\[6pt]
      ZFP, $\eta = 10^{-3}$ &  1.05  & 1.16  &  0.10  & 11.2 \\[6pt]
      ZFP, $\eta = 10^{-2}$ &  1.06  & 1.17  &  0.12  & 17.6 \\[6pt]
    \end{tabular}
    \vspace*{0.2cm}
    \caption[Earth domain results]{
      Comparison of ZFP compression to autoencoder compression for various
      compression tolerances. There is not a significant difference in the
      speedup of the autoencoder compression vs the ZFP compression,
      but the autoencoder is able to achieve a much higher compression ratio (and hence much smaller memory footprint)
      while still maintaining solution accuracy. Relative $l_2$ error is
      computed with respect to the checkpoint solution. 
    }
    \label{table:zfp_vs_ae}
\end{table}

\subsection{DIAS regularization with autoencoder compression}
As a capstone result, we show how our proposed autoencoder compression
approach can be combined with the DIAS regularization \cite{nguyen2022dias} approach
on the earth-scale seismic inverse problem. 
While the DIAS algorithm has the benefit of only applying regularization
in the inactive subspace, speed is not one of it's strengths.
Simply estimating the active subspace for a nonlinear inverse problem
requires numerous samples of the gradient --- each evaluation of which
requires 2 PDE solves. In the case of large-scale seismic inversion,
each gradient evaluation can take several minutes on dozens of compute nodes.
We propose to combine the DIAS regularization algorithm with
autoencoder compression
to alleviate the high cost and memory footprint of checkpointing in estimating the active subspace. Since we have already discussed the speedup and memory footprint  gained by the autoencoder compression at length above, we will focus on the solution quality using DIAS plus autoencoder compression in the following.

Let us briefly recap the DIAS regularization. We will focus on the DIAS-F (see \cite[Section 4]{nguyen2022dias} 
variant, which uses the full data misfit and only applies regularization in
the inactive subspace.
First, the active subspace is estimated via,
\begin{equation}
  \begin{bmatrix}
    \Wbo &\Wbt
  \end{bmatrix}
  \begin{bmatrix}
    \Lambda_1 & 0 \\
    0 &\Lambda_2
  \end{bmatrix}
  \begin{bmatrix}
    \Wbo^T \\ \Wbt^T
  \end{bmatrix}
  = \Cb := \int \grad_{\qoi} f\LRp{\qoi} \grad_{\qoi} f\LRp{\qoi}^T \probMeas[\text{prior}][\qoi] d \qoi,
\end{equation}
where $f:= \half \norm{\pto\LRp{\qoi} - \obs}^2_{\invNoiseCov}$, $\Wbo$ are the eigenvectors of
$\Cb$ corresponding to the active subspace, and $\Wbt$ are the eigenvectors corresponding to the
inactive subspace. As in Section \ref{section:training_data}, 
$\probMeas[prior][\qoi]$ is the prior probability density from the Bayesian formulation
of the inverse problem. 
The DIAS solution is then defined to be
\begin{equation}
  \label{eqn:dias_prior_opt}
  \qoiDias = \underset{\qoi}{\text{argmin }} \half \norm{\pto \LRp{\qoi} - \obs}^2_{\noiseCov^{-1}}
  + \half \norm{\Wbt^T \LRp{\qoi  - \qoi_0}}^2_{\LRp{\Wbt^T \priorCov \Wbt}^{-1}},
\end{equation}
that is, it is a Tikhonov solution with regularization acting only on the inactive space.

For large scale problems, it is infeasible to compute and store the large matrix $\Wbt$, since the
active subspace tends to be much smaller than the inactive subspace.
To avoid explicitly computing and storing the potentially large matrix $\Wbt$, it can be shown
using \cite[Theorem 5.8]{schott2016matrix} that
\begin{equation}
  \norm{\Wbt^T \LRp{\qoi  - \qoi_0}}^2_{\LRp{\Wbt^T \priorCov \Wbt}^{-1}} =
  \norm{\proj \LRp{\qoi  - \qoi_0}}^2_{\LRp{\proj \priorCov \proj}^{\dagger}}
\end{equation}
where $\proj := \Wbt \Wbt^T = \ident - \Wbo \Wbo^T$ is a projection matrix onto the inactive subspace
and $\LRp{\proj \priorCov \proj}^{\dagger}$ denotes the pseudoinverse.

Lest we conclude that we have a free lunch in computing the regularization term arising
from the active subspace prior via $\Wb_1$, there is one more challenge that remains.
Note that we need to apply not the prior covariance, but its (pseudo-)inverse.
That is, we need to compute $\LRp{\proj \priorCov \proj}^{\dagger}$.
There are several factors at play here.

\begin{enumerate}
\item $\priorCov$ is often a sparse matrix, while $\proj \priorCov \proj$ is very likely a dense matrix.
\item We need to invert the dense matrix $\proj \priorCov \proj$ using the pseudo-inverse.
\item This matrix is in the dimension of the parameter,  $n$. 
\end{enumerate}

For cases where it is feasible to form and invert $\proj \priorCov \proj$, it is also
feasible to simply compute $\Wb_2$ explicitly. Then we can work with the lower
dimensional form given in \eqref{eqn:dias_prior_opt} which still involves inversion,
but of a lower dimensional, full rank matrix. For cases where this is not feasible,
we need to develop more sophisticated techniques for computing or estimating
$\proj \priorCov \proj$.

The first thing one might try is to find an identity that allows us to break the
projection out of the pseudo-inverse, leaving $\priorCov$ by itself.
We might hope that
\begin{equation*}
  \LRp{\proj\priorCov \proj}^{\dagger} = \proj \priorCov^{-1} \proj,  
\end{equation*}
but this is unfortunately not the case, as we show using Schur complements. 
Decompose the prior covariance as
\begin{align*}
  \priorCov = 
  \LRs{\begin{array}{c|c}
         \Wb_1^T \priorCov \Wb_1 & \Wb_1^T \priorCov \Wb_2 \\
         \hline
         \Wb_2^T \priorCov \Wb_1 & \Wb_2^T \priorCov \Wb_2
      \end{array}}
\end{align*}
and the inverse prior covariance as
\begin{align*}
  \priorCov^{-1}= 
  \LRs{\begin{array}{c|c}
         \Wb_1^T \priorCov^{-1} \Wb_1 & \Wb_1^T \priorCov^{-1} \Wb_2 \\
         \hline
         \Wb_2^T \priorCov^{-1} \Wb_1 & \Wb_2^T \priorCov^{-1} \Wb_2
      \end{array}}.
\end{align*}
Using Schur complements gives
\newcommand{\Amat}{\Wb_1^T \priorCov \Wb_1}
\newcommand{\Bmat}{\Wb_1^T \priorCov \Wb_2}
\newcommand{\Cmat}{\Wb_2^T \priorCov \Wb_1}
\newcommand{\Dmat}{\Wb_2^T \priorCov \Wb_2}
\begin{align}
  \label{eqn:cov_approx_not_equal}
  \Wb_2^T \priorCov^{-1} \Wb_2 &= 
  \LRs{\Dmat - \Cmat \LRp{\Amat}^{-1} \Bmat}^{-1}.                               
\end{align}
While the inverse can be expanded using the Sherman-Morrison-Woodbury formula,
the important observation is that the second term,
\begin{equation*}
  \Cmat \LRp{\Amat}^{-1} \Bmat
\end{equation*}
is not necessarily 0, except when $\priorCov$ is scaled identity.
Thus, making the approximation $\LRp{\proj\priorCov \proj}^{\dagger} = \proj \priorCov^{-1} \proj$
incurs some error in the general case.
Even though there is some error, it becomes necessary from a practical perspective
to use the approximate form $\proj \priorCov^{-1} \proj$ for large scale problems. 

Before the DIAS algorithm can be applied to the seismic inverse problem,
we must first define an algorithm for applying DIAS to nonlinear inverse problems.
One advantage of using the active subspace to define the data-informed
subspace is that the active subspace is formed from the \textit{global} average of the
outer product of the gradient. This would indicate that the DIAS algorithm
can be applied naively to nonlinear inverse problems:
\begin{equation}
  \qoiDias = \underset{\qoi}{\text{argmin }} \half \norm{\pto \LRp{\qoi} - \obs}^2_{\noiseCov^{-1}}
  + \half \norm{\proj \LRp{\qoi  - \qoi_0}}^2_{\LRp{\proj \priorCov \proj}^{-1}}.
\end{equation}
However, the fact that the gradient samples are computed from draws of the prior
biases the active subspace toward the directions in which the misfit function $f$
is most sensitive near the prior mean, $\qoi_0$.

Recalling that the goal of the DIAS algorithm is to remove regularization in directions that the data
are more informative and that the DIAS solutions are likely closer to the usual Bayesian solution
than to the prior mean for a well-chosen prior,
we propose a two-step algorithm:
\begin{enumerate}
\item Approximately solve the usual MAP problem for $\qoi_{\text{MAP}}$ using the full prior:
  \begin{equation}
    \qoi_{\text{MAP}} = \underset{\qoi}{\text{argmin }}
    \half \norm{\pto \LRp{\qoi} - \obs}^2_{\invNoiseCov}
    + \half \norm{\qoi - \qoi_0}^2_{\invPriorCov}.
    \label{eqn:MAP_solution_dias}
  \end{equation}
\item Compute the active subspace centered at $\qoi_{\text{MAP}}$ and, with initial guess $\qoi_{\text{MAP}}$, solve DIAS problem:
  \begin{equation*}
    \qoiDias = \underset{\qoi}{\text{argmin }}
    \half \norm{\pto \LRp{\qoi} - \obs}^2_{\noiseCov^{-1}}
    + \half \norm{\proj \LRp{\qoi  - \qoi_0}}^2_{\proj \invPriorCov \proj}.
  \end{equation*}
\end{enumerate}

\begin{figure}
  \centering
  \begin{subfigure}{0.3\textwidth}
    \includegraphics[trim=4.7cm 4.7cm 4.7cm 4.7cm,clip=true,width=\textwidth]{figures/uqearth/ae_solution_true_scale}
    \caption{MAP}
  \end{subfigure}
  \begin{subfigure}{0.3\textwidth}
    \includegraphics[trim=4.7cm 4.7cm 4.7cm 4.7cm,clip=true,width=\textwidth]{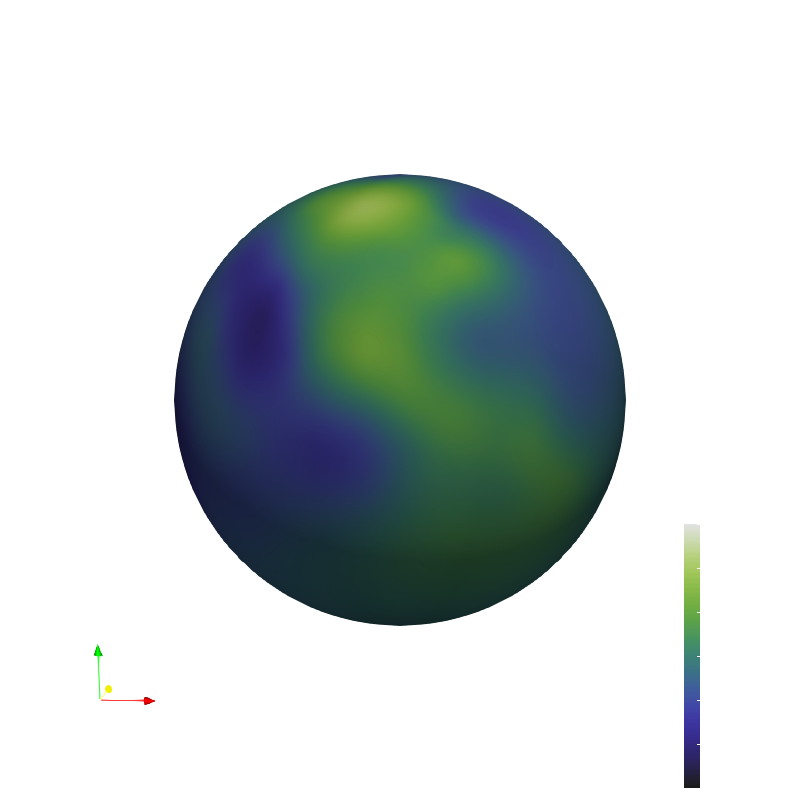}
    \caption{DIAS}
  \end{subfigure}
  \begin{subfigure}{0.3\textwidth}
    \includegraphics[trim=4.7cm 4.7cm 4.7cm 4.7cm,clip=true,width=\textwidth]{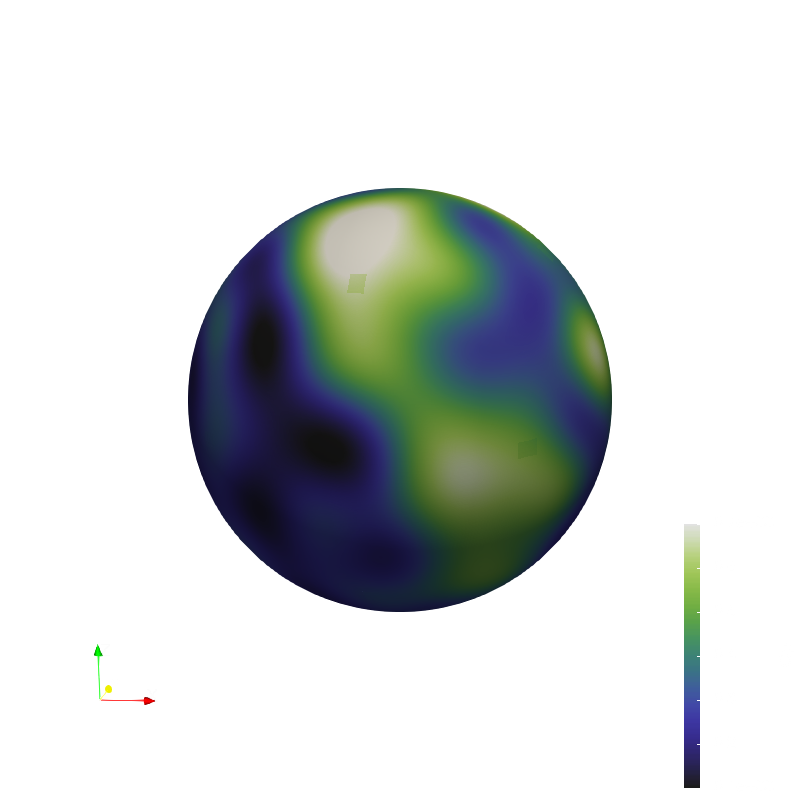}
    \caption{True}
  \end{subfigure}
  \begin{subfigure}{0.085\textwidth}
    \begin{tabular}{c}
        \hspace*{-0.2cm}
        \includegraphics[trim=33.2cm 0cm 1.51cm 12.5cm,clip=true,width=\textwidth]{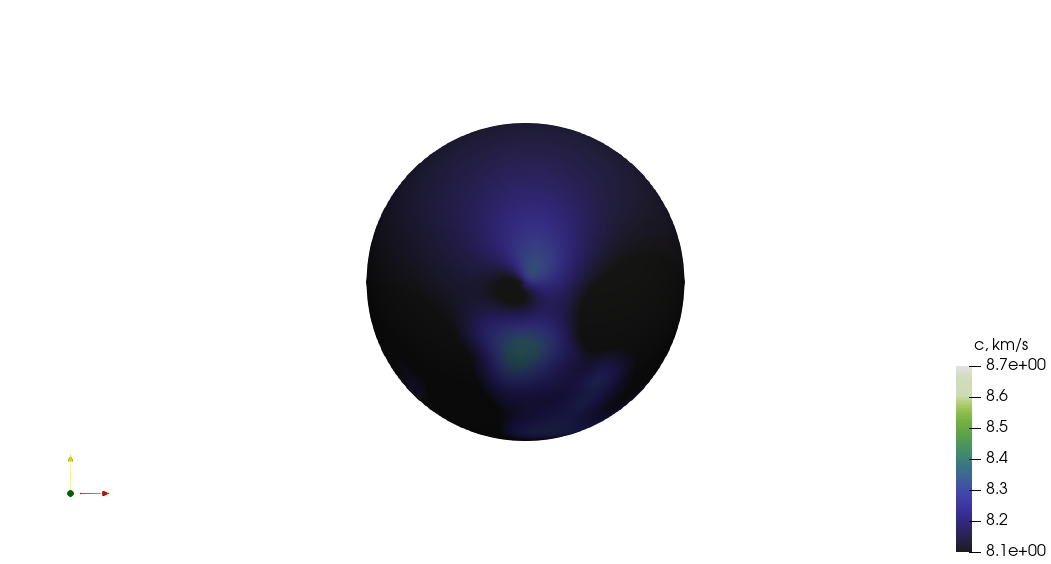}\\
        $c, $ km/s
    \end{tabular}
  \end{subfigure}
  \caption[MAP vs DIAS earth solutions]{
    The DIAS solution does not differ noticeably from the MAP solution
    for the seismic inverse problem. Results shown for mesh
    with 11,166,784 DoFs. 
  }
  \label{fig:dias_earth}
\end{figure}

Even with acceleration via autoencoder compression,
the cost of gradient evaluations is still high so we
estimate the active subspace using 30 gradient samples and take the
active subspace to have dimension 5.
The difference between the DIAS and the MAP estimate are small in this
case. Compared to the true solution, the MAP estimate has average relative error of
0.3293\% while the DIAS refined solution has average relative error of
0.3291\%. This is reflected visually in Figure \ref{fig:dias_earth}.

\section{Conclusions}
In this paper, we proposed an autoencoder compression algorithm to mitigate the high
storage requirements for solving large-scale time-dependent PDE constrained
inverse problems. We show the feasibility and scalability of this approach against the popular checkpointing strategy for a
seismic inverse problem on both a simple box domain with uniform refinement
and a complex spherical domain with adaptive mesh refinement and a nonconforming mesh.
In order to effectively train the autoencoder, we present a
data generation procedure based on the Bayesian formulation
that provides a problem-focused data generation scheme. We proposed two separate autoencoder compression algorithms for spatial and temporal domains.
The scalability of each algorithm was shown,
providing a clear close-to-ideal speed advantage on large scale problems compared to the
traditional checkpointing approach and a substantial improvement in the
compression ratio of the state-of-the-art floating point compression
package, ZFP. While there are clear advantages of using autoencoders
to perform compression, care must be taken in their implementation
in order for there to be an advantage over the checkpointing approach. As an important application of the proposed autoencoder algorithm, 
we combined the proposed autoencoder compression approach
with the DIAS prior to show how the DIAS method
can be affordably extended to large scale problems without the need of checkpointing and large memory.


\ifx\anonymousAuthor\undefined
\section{Acknowledgements}
We would like to thank Georg Stadler for his many correspondences and helpful 
advice on stably solving the forward and inverse problem in large scale runs;
Jau-Uei Chen for sharing his insight on the discontinuous Galerkin
method; Sheroze Sheriffdeen and Hwan Goh for their many discussions on
autoencoders, compression, and HPC;
and the Texas Advanced Computing Center (TACC) at the 
University of Texas at Austin for providing 
HPC resources that contributed to the results presented in this work. 
URL: http://www.tacc.utexas.edu. 
This research is partially funded by the National Science Foundation awards
NSF-OAC-2212442, NSF-2108320, NSF-1808576 and NSF-CAREER-1845799; 
and by the Department of Energy awards DE-SC0018147 and DE-SC0022211.
\fi

\bibliography{references}

\end{document}